\documentclass[a4paper,12pt]{article}
\usepackage{localrus}
\title{Израиль Моисеевич Гельфанд}
\author{Александр Шень\thanks{Alexander Shen, LIRMM, Univ Montpellier, CNRS, Montpellier, France,\\ \protect\url{email: sasha.shen@gmail.com}}}
\date{}
\begin{document}
\maketitle

\begin{abstract}
Я имел возможность (пусть немного) встречаться и даже сотрудничать с Израилем Моисеевичем Гельфандом. В этом тексте (начатом в 2003--2013 годах) я пытался вспомнить некоторые свои впечатления об этих встречах и сформулировать некоторые свои соображения и наблюдения по части преподавания математики, в значительной степени сложившиеся под его влиянием (хотя, конечно, И.\,М. за них ответственности не несёт).
\end{abstract}

У Израиля Моисеевича было много учеников и сотрудников, которые знали его много дольше и много лучше, чем я, и могли бы рассказать больше и интереснее --- но пока далеко не все это сделали. Имея это в виду, я всё же попытаюсь вспомнить о своём знакомстве с И.\,М. и заранее прошу прощения, если что-то перепутал или недопонял.

Про легендарный семинар Гельфанда (математический) я слышал рассказы от своих знакомых, но, конечно, никогда бы не решился туда пойти --- по недостатку квалификации и образования. Однако в середине 1980-х годов Алексей Львович Семёнов и Анатолий Георгиевич Кушниренко позвали меня участвовать в написании учебников по информатике, и в это же время Гельфанд решил устроить небольшой семинар по темам, связанным с компьютерами, с участием Кушниренко. Одновременно Гельфанд пытался как-то освежить и наладить работу ВЗМШ (Всесоюзной заочной математической школы, которую он организовал в середине 1960-х), а я там многих знал и имел некоторый опыт работы в качестве студента-проверяющего (а кроме того, преподавал в математических классах в 91 и 57 школах), и это было другой причиной знакомства.

\subsection*{Семинар по информатике}

Семинар Гельфанда по <<информатике>> (не помню, употреблял ли это слово тогда сам И.\,М.) был небольшим. Он проходил в лабораторном корпусе <<А>> МГУ, и скорее был посвящён знакомству с самыми разными работами, чем какой-то конкретной области. Как я понимаю, до этого Гельфанд с сотрудниками много занимался применением компьютеров в медицине и вообще тем, что называется <<распознавание образов>>. В основном Гельфанда в информатике интересовали не математические задачи (хорошо поставленные, где можно доказывать какие-то теоремы, например, нижние или верхние оценки времени работы алгоритмов и так далее --- эта наука и тогда уже была довольно развита, а с тех пор там много новых замечательных продвижений), а именно практические, нечёткие, плохо поставленные задачи, и математические методы их решения. Например. обсуждался <<метод оврагов>> поиска минимума функций большого числа переменных: если график функции имеет вид <<оврага>>, то не годится ни большой шаг в направлении градиента, который выведет далеко за пределы оврага, ни маленький, который не позволит дойти по оврагу до минимума --- надо сочетать то и другое. Помню ещё, как нам рассказывали о работах группы Шели Айзиковича Губермана по распознаванию рукописных букв.\footnote{%
     Потом я ходил смотреть на эти программы --- уже без Гельфанда --- в <<компьютерный клуб>>, который был создан по инициативе Гарри Каспарова (и при поддержке Евгения Павловича Велихова) группой энтузиастов, в которую входил Степан Пачиков (будущий председатель Совета клуба), Семёнов, Александр Звонкин (автор <<Дневника математического кружка>>), Андрей Тоом (<<алгоритм Тоома--Кука>>, <<правило Тоома>> коррекции ошибок для клеточных автоматов) и многие другие. У меня нашёлся [написано в 2007 году --- вероятно, этот список сохранился и теперь, но остался в Москве] список телефонов от 23.11.1987; в нём $43$~<<члена>> клуба, $5$~<<кандидатов в члены>>, $3$~<<члена-корреспондента>>, включая меня самого, и $7$~<<сочувствующих и присматривающихся>>.
     
      Клуб этот сначала был на Арбате, потом (с помощью Семёнова) был получен особняк на Рождественском бульваре, правда, в ужасном состоянии и с невыехавшим общежитием какого-то медицинского училища. Один раз, когда туда должно было прийти какое-то начальство, были закуплены плакаты о научно-техническом прогрессе для завешивания дырок в стенах. Компьютеры (фирмы Atari) для этого клуба привёз Гарри Каспаров (за счёт своих шахматных гонораров, видимо), и туда приходили заниматься (не просто играть, а именно заниматься) школьники разного возраста; одним них был Миша Белкин --- потом ученик 57 школы, а ныне известный специалист по Computer Science (\url{https://scholar.google.com/citations?user=Iwd9DdkAAAAJ&hl=en}). 
     
Распознавание рукописных букв, к моему удивлению, происходило довольно хорошо. Потом из работ этой группы получились первые изделия организованной Пачиковым (с товарищами) фирмы <<Параграф>>, которыми заинтересовалась фирма Apple.}
Ещё на семинаре рассказывали о дифференциальной диагностике болезней лёгких (распознавание опухолей по рентгенограммам) и о многом другом.

Конечно, в сравнении с математическом семинаром (о котором речь пойдёт дальше) роль этого  <<информатического>> семинара была совсем иной. Сложность с плохо поставленными задачами в том, что очень трудно передавать накопленный опыт, и потому каждый следующий раз приходится начинать почти с того же места. Кроме того, ими нельзя заниматься теоретически, нужен реальный проект с реальными заказчиками.\footnote{%
     Может быть, в наше время, когда с появлением поисковых машин и прочего интернета спрос сильно возрос, продвижение ускорится. [Написано в 2007 --- сейчас прогресс ускорился настолько, что вообще происходит непонятно что, в частности, программы делают то, что они явно делать не могут, и т.\,п.]} 
     Так что сам семинар был скорее местом знакомства с интересными людьми, чем способом совместной работы его участников.\footnote{И.\,М. много работал со своими сотрудниками по части распознавания образов и применения компьютеров --- в одном из препринтов этой группы, насколько я помню, были даже разные <<пословицы и поговорки>>, возникшие в ходе этой работы. Помню одну из них: <<Не то диво, как мужик на луну залез, диво, как он с неё слезать будет. Записать на ленту --- не фокус, фокус с ленты прочитать>>. Но эта работа происходила до/вне семинара.}

\subsection*{Заочная школа и книжка по алгебре}

В Заочной школе Гельфанд хотел освежить методические материалы.
В самом начале работы ВЗМШ с его участием были написаны книжки
(<<Метод координат>>, <<Функции и графики>>, соавторами
Гельфанда были Елена Георгиевна Глаголева, Эммануил Эльевич
Шноль и Александр Александрович Кириллов-старший), которые
много лет использовались (и продолжают использоваться), но было
видно, что этого мало. Одна из книжек, начатых по этому поводу
(<<Тригонометрия>>), была потом дописана уже без участия И.\,М.
его соавторами Сергеем Михайловичем Львовским и Андреем
Леоновичем Тоомом (текст см. в~\url{https://www.mccme.ru/free-books/lvovski/trig.pdf}). А сам И.\,М. уже в Америке в 2001 году выпустил книжку по тригонометрии с Марком Солом (ISBN 978-0817639143). (Обсуждение содержания этих книжек и ссылки см. ниже.)

 В другой книжке (<<Алгебра>>) Гельфанд
предложил участвовать мне, и через много лет она-таки вышла, хотя и позже всех сроков --- она была объявлена в планах издательства то
ли на 1988, то ли на 1989 год, а вышла только в середине 1990-х
годов в английском варианте в Америке, русское издание (текст см. в~\url{https://hal-lirmm.ccsd.cnrs.fr/lirmm-01486516/document}) было ещё
позже.

Работа над книжкой начиналась так: я приходил домой к Гельфанду
(он жил около Юго-Западной), и И.\,М. объяснял мне, чт\'{о} нужно
написать, почти диктуя текст. Записав и отредактировав этот
текст, я приходил к нему снова, он смотрел текст и делал разные
замечания, после чего текст переписывался заново (часто совсем
по-другому), и так повторялось несколько раз. Одновременно
И.\,М. пытался объяснить мне, чего бы хотелось и как надо писать,
и постепенно я научился лучше понимать, что ему нравится и что
не нравится, так что количество итераций уменьшалось. Вторую
половину книжки я уже писал после общего обсуждения содержания,
и он только делал замечания по готовому тексту. В Москве так это
до конца и не дошло, и будучи в 1991 году в Бостоне (в MIT), я
ездил к И.\,М. (в Rutgers University, видимо) для продолжения
работ.
\smallskip

Помню разные замечания ИМ из этих разговоров:

\begin{itemize}

\item Написав абзац, спросите себя: что мы хотим сообщить
читателю в нём? если ответ неясен, абзац можно вычеркнуть. Если
ответ ясен, то абзац тоже можно вычеркнуть, заменив на этот
ответ.

\item Книжки можно писать по-разному, и эту книжку надо писать
скорее подражая Пушкину, чем Толстому.

\item Говоря об попытках преждевременных объяснений в
преподавании, И.\,М. любил рассказывать такой анекдот. Ребёнок
спрашивает: <<Мама, что значит \emph{аборт}?>> --- <<Ну как тебе
сказать\ldots>> --- смущается мама и начинает что-то сбивчиво
объяснять. Ребёнок слушает, слушает и наконец не выдерживает:
<<При чём тут это?! я тебя спрашиваю, что значит: \emph{волны бьются
о борт корабля}?>>

\item Барток в <<Микрокосмосе>> создал цикл пьес для детей,
каждая из которых иллюстрирует какой-то приём в миниатюре. Так
же надо делать и в книжках по математике для школьников. (Вообще
И.\,М. часто приводил музыкальные примеры --- помню, как он восхищался
Моцартом и записью <<Дон Жуана>> с Фуртвенглером, и этот
пример был мне понятнее, чем Барток, которого я почти не знал.)

\end{itemize}

В то время дочка Гельфанда Таня была как раз в младших классах,
и его разговоры с ней повлияли на содержание книжки (и даже
упоминаются в ней).

В середине 90-х годов Гельфанд организовал в Ратгерсе некий американский вариант заочной школы (под названием Gelfand Outreach Program in Mathematics).  В качестве заданий использовались \emph{Метод координат} и \emph{Функции и графики}, переведённые раньше, и туда же предназначался английский перевод книжки по алгебре.
Кроме того, он хотел написать книжку по
геометрии, и какие-то материалы были написаны им с женой
(Татьяной Алексеевской), но это было только самым началом
работы --- их книжка вышла через много лет, в 2020 (ISBN 978-1071602973). Наконец, требовалось составлять отдельные задания для
этой школы, и я в этом
тоже пытался помогать, приехав в Ратгерс на несколько месяцев по его
приглашению. (Непосредственно в работе
школы, то есть в проверке заданий, я не участвовал.)

Книжка по геометрии с моим участием так и не была написана; было
много набросков, но ни одним из них И.\,М. не был доволен. По его
предложению я стал подбирать и группировать задачи и снабжать их
комментариями, и это ему понравилось больше. Я продолжал уже в
Москве, и окончательного одобрения И.\,М. текст так не получил.
После разных доделок и переделок он был опубликован в 2013 году (текст см. в~\url{https://hal-lirmm.ccsd.cnrs.fr/lirmm-01235075v1/document}).

Параллельно И.\,М. предложил мне некоторую тему для математических
занятий (элементарная некоммутативная алгебра, свободные
некоммутативные поля), но из этого ничего не вышло (прежде всего
по моей неспособности).

\subsection*{Математический семинар}

Ещё в Москве Гельфанд позвал меня на знаменитый математический
семинар (по понедельникам в аудитории 14-08 главного здания МГУ),
и я имел возможность наблюдать за тем, как он происходит. Как
уже многие писали, он совсем не был похож на стандартный
(особенно западный) семинар, где докладчик заранее готовит
доклад (часто в виде слайдов, теперь компьютерных) и потом его
воспроизводит.

Трудно было сказать заранее, когда начнётся доклад (до него были
разные обсуждения), не всегда было ясно, кто докладчик. И уж
никто точно не знал, как пойдёт доклад, сколько он продлится и
чем кончится. Для непривычных докладчиков это, конечно, было
некоторым шоком, но смысл в этом был: ведь обычно докладчик
лучше знает предмет доклада, чем слушатели, и это ему скорее
следует приспосабливаться к слушателям, чем наоборот. С другой
стороны, неожиданный для докладчика взгляд на предмет доклада со
стороны продвинутых слушателей может быть интересным и для него
самого.

На семинаре были люди с разной подготовкой, и это не мешало им
извлекать пользу --- в частности, из <<лирических
отступлений>> И.\,М., который часто по ходу дела комментировал
разные понятия и факты, имеющие (или даже не очень имеющие)
отношение к теме доклада. И.\,М. вполне мог спросить кого-то из
слушателей (часто начинающих), понимают ли они происходящее, и
даже вызвать их к доске пересказывать (если поняли) --- в качестве
<<контрольных слушателей>>. Как-то раз я попал в их число,
и Максим Концевич (ученик Гельфанда и постоянный участник
семинара) объяснил мне, что не надо удивляться --- <<ты сейчас
исполняешь роль первокурсника>>.

Кстати, про выступление М.К. на семинаре Гельфанда рассказывали
байку --- якобы И.\,М. несколько раз просил его говорить громче, и
наконец Максим напрягся и прокричал: <<Я громче не могу!>> ---
<<Что? Не слышу!>> --- отвечал И.\,М. Но при этом я сам не был.
Зато я застал другой случай, когда после выступления Максима
И.\,М. сказал участникам что-то вроде <<теперь можете его
спросить, пользуйтесь случаем --- через несколько лет он будет
знаменитым, и это будет не так просто>> (тогда он был
аспирантом). Другой раз, когда я по какому-то поводу сказал
Гельфанду в разговоре, что мало что понимаю, и это особенно
заметно в сравнении с Сашей Полищуком и Лёней Посицельским
(которых я знал ещё школьниками, а к моменту разговора они были
студентами мехмата), И.\,М. улыбнулся и ответил: <<Нашли с кем
сравнивать>>.

\subsection*{Сотрудники и коллеги}

И.\,М. внимательно относился к людям и ценил их достоинства, но не
церемонился с ними, и это бывало источником обид. Как-то я
присутствовал при такого рода разговоре (уже в Америке,
кажется --- не помню, кто был собеседником), и когда обиженный
собеседник ушёл, И.\,М. сказал мне примерно так: <<бывают такие
обидчивые люди, что с ними невозможно иметь дело; единственный
способ --- не церемониться с ними, и если они привыкнут, то уже
можно и работать>>. Но способ этот годился не всегда, и весьма
достойные люди через много лет и в другой части света оставались
в обиде на И.\,М. (например, А.Л.\,Тоом).

Помню ещё, что когда И.\,М. поручал мне координировать какие-то
работы для заочной школы (в Москве), он сказал: <<поторопите
их и скажите, что завтра я вам позвоню и буду спрашивать, что
сделано>> --- а после некоторого моего замешательства добавил:
<<хотите, действительно вам позвоню?>>.

Ещё он как-то сказал, что хорошо, когда у человека больше
обязанностей, чем физических возможностей и времени --- тогда он
занимается тем, что ему действительно интересно.

Можно успеть больше, говорил И.\,М., если вместо отдыха
переключаться на другую деятельность (как от математики к
биологии). Говоря об отъезде в Rutgers University, И.\,М. сказал, что это
добавило ему несколько лет активной жизни. В середине 1990-х
работа Гельфанда там выглядела так: в его кабинет приходили
коллеги (математики и биологи) и начинались обсуждения
сделанного ими после предыдущего разговора; никаких записей И.\,М.
не делал; часто чередовались самые разные темы, и можно было
только удивляться, как И.\,М. всё это понимает\ldots

Не заботясь о самолюбии своих коллег и сотрудников, И.\,М. часто и
очень существенно им помогал. Когда случилась тяжелейшая авария
с Юрой Кушниренко (сыном А.Г.\,Кушниренко), Гельфанд решительно
вмешался в ситуацию, поднял по тревоге своих медицинских
знакомых и, видимо, в значительной мере спас ему жизнь. Он
помог Марату Ровинскому перевестись на мехмат (договорившись с
Садовничим, вероятно). Думаю, что люди, хорошо знавшие И.\,М.,
могут вспомнить сотни (если не тысячи) таких случаев.

\subsection*{Анекдоты от Гельфанда}

Говорят, что Кириллов (старший) собирался выпустить сборник
анекдотов (в старинном и современном смысле слова) <<от
Гельфанда>>. Кажется, пока такого сборника нет, так что
попытаюсь пересказать некоторые.

\begin{itemize}

\item Летит самолёт, в нём волк, медведь и ворона. Ворона
крутит штурвал, изображая пилота. Волк её спрашивает: <<Ты что
делаешь?>> --- <<Выпендриваюсь>>. Волк: <<Дай-ка я
повыпендриваюсь>>. (Самолёт немного трясёт.) Наконец, садится
медведь, самолёт болтается туда-сюда и наконец разваливается.
Ворона летит и думает: <<Странное дело, летать не умеют --- а
выпендриваются>>.

\item Когда в вычислительное бюро, которым заведовал Меир
Феликсович Бокштейн (имени которого гомоморфизм), завезли трофейные
немецкие счётные машинки, М.\,Ф. первым делом решил попробовать
разделить на нуль. Машинка щёлкала, щёлкала, и наконец каретка
вылетела. С тех пор стало ясно --- на нуль делить нельзя!

\item Медведь --- лисице: <<Приходи ко мне, я тебя съем!>> --- <<А когда?>> --- Медведь, делая пометки в записной книжке:  <<Сегодня вечером>>. (Лисица в ужасе уходит.) Медведь встречает волка: <<Приходи ко мне, я тебя съем!>> (Глядит в книжку.)
<<Сегодня придёт лисица, так что ты приходи завтра утром>>. (Волк
в ужасе уходит.) Медведь видит зайца: <<Приходи, я тебя
съем!>> (Делает пометки.) <<Сегодня и завтра утром лисица и
волк, так что приходи завтра на обед.>> Заяц: <<А пошёл ты
$\langle\ldots\rangle$!>> --- <<Ну ладно, не хочешь ---
вычёркиваю.>>

\item Мадам --- привередливому клиенту: <<Ну прямо и не знаю,
кто вас устроит, вы от всех отказываетесь. Могу разве что
предложить свои услуги>>. --- <<Но на вас, мадам, целое
заведение\ldots>> --- <<Иногда так устаёшь от
оргработы\ldots>>.

\end{itemize}

\subsection*{Уроки Гельфанда}

Название этого раздела не совсем соответствует его содержанию --- когда я стал вспоминать, чему и как я научился у И.\,М., то понял, что эти наблюдения и выводы трудно разделить с моими собственными соображениями --- так что не буду их разделять и сразу скажу, что И.\,М. за них ответственности не несёт. 

Традиционно (и сейчас тоже часто так бывает) лекции и научные семинары по математике проходили как монолог лектора у доски или у подготовленной им презентации\footnote{Часто докладчики предпочитают слайды, потому что иначе они не успеют написать всё нужное на доске --- не учитывая, что важно не только написать, но и прочитать. Попытки что-то дополнительно разъяснить, пользуясь слайдами, часто выглядят как сбор грибов с дрезины: можно ездить туда-сюда, но нельзя сойти с рельсов и достать гриб.} со слайдами (когда-то на <<прозрачках>>, теперь компьютерными) --- а слушатели, даже если их и призывали задавать вопросы по ходу дела (а не в конце), в основном молчали. Сейчас, когда можно найти и прочесть статью, о которой рассказывается в докладе, или любой из десятков учебников, или даже видео каких-то лекций и семинаров на ту же тему, или, наконец, видео прямо этого же выступления (потом), смысл такой церемонии становится менее ясным. 

Что слушатель получает дополнительно, придя на лекцию? (Что теряет --- понятно: он должен сидеть и слушать в темпе рассказчика, не имея возможности пропустить неинтересное место, или обдумать непонятное, или поставить видео на скорость $1.5$ или $0.75$, или просто налить себе чаю.)  А что получает докладчик? Иногда такое выступление воспринимается как <<апробация>> результатов докладчика --- успешная, если все разошлись в почтительном молчании или задали какие-то <<вопросы для приличия>>. (А если кто-то начал что-то уточнять, не понявши, или говорить, что ещё можно или нужно делать и как --- то это личный враг, срывающий апробацию.) Но и это не очень понятно зачем, если есть \texttt{arxiv.org}.

Преимущества возникают, если рассказывающий ориентируется на реакцию слушателей и понимает, что важно не то, что он расскажет, а то, что они поймут. Лектору в этом отношении приходится рассчитывать только на самого себя, а на семинаре есть руководитель, который обычно лучше знает участников и может понять их реакцию. Даже простой вопрос <<кто понял>> с просьбой к понявшим поднять руки и при необходимости пересказать другим, может вернуть докладчика к  (прискорбной) реальности\footnote{Можно учиться рассказывать методом проб и ошибок --- но для этого надо уметь увидеть ошибку, что трудно сделать, рассказывая поверх презентации.} --- и, например, побудить его рассказывать совсем не то, что он собирался, а то, что лучше для слушателей. Иногда Гельфанд просил вместо общих определений разобрать какие-то простые примеры. Или наоборот --- когда докладчик долго пытался что-то пояснить неформально, И.\,М. мог сказать, что <<невозможно мотивировать определение, которое не дано>>. 

Конечно, сотрудничество требует смирения не только от докладчика, но и от слушателей --- если они стараются показать, какие они сообразительные, вместо того чтобы сказать, что ничего не поняли, когда не поняли, и промолчать, когда поняли, то и самый лучший докладчик будет в тупике. (Руководители семинара могут помочь, задавая полезные вопросы и препятствуя вредным.) 

Наконец, надо учитывать, что слушатели разные --- нужно, чтобы те, кто может многое понять, имели такую возможность, так что нельзя ориентироваться только на минимальный уровень. (И, конечно, хуже всего, если все уже всё поняли, а руководители ещё чего-то хотят от докладчика --- так тоже бывает.)  Про лекции для студентов И.\,М. говорил, что там должно быть что-то понятное и интересное для студентов разных уровней. При семинаре И.\,М. иногда бывали параллельные <<семинары для начинающих>>, где более продвинутые слушатели делились своими знаниями с менее продвинутыми. И, пожалуй, отношение И.\,М. к начинающим было гораздо более деликатным, чем к <<маститым учёным>> --- пример, которому многим бы стоило следовать в его позитивной части.

Мой товарищ Андрей Коган, который в какой-то момент приходил устраиваться (тогда говорили именно так, а не <<наниматься>>) на работу в группу Гельфанда и Губермана, так описывал в 2003 году (\url{https://kogan.livejournal.com/4724.html}) свою встречу с Гельфандом (Коган в момент встречи был недавним выпускником физтеха):

\begin{quote}
Народ тусовался перед аудиторией, но почему-то не входил. В глаза бросался крайне солидного вида крупный мужчина в отлично выглаженном костюме. Ясно было, что он тут главный. Подумав, я решил, что Гельфандом (коего я раньше не видел) этот солидный человек все же не является. Тут как раз в коридоре материализовался одуванчикового вида старичок в сопровождении небольшой свиты. Я представил, как он сейчас со светлым и чуть рассеянным взором пройдет сквозь нас в аудиторию. Но промазал, ибо старичок повернулся к тому самому солидному пиджаку, ткнул в него пальцем и произнес необычайно сварливо и безо всяких знаков препинания:

--- А почем\'{у} Вы не пригласили на сегодняшний семинар такого-то-такого-то --- чт\'{о} ?!

--- Я, Израиль Моисе\ldots

--- Помолчите!

Через пять минут, доведя солидного господина до состояния полного фарша, он обратил, наконец, внимание на нас (не помню, кто меня привел, кажется, Дзюба [один из сотрудников И.\,М.])

--- Пойдемте!

Кажется, было что-то вроде кабинета.

--- Чем Вы занимаетесь? (Имелась в виду дипломная работа.)

Я рассказал (не помню, о чем, --- честно говоря, я ничем особенным тогда и не занимался).

--- А если мы Вам скажем, что это чушь?

Тут я заметно обиделся и процедил что-то вроде <<я подумаю над Вашими словами>>. Впрочем, долго мы на эту тему не беседовали и перешли к математике.

Любопытно, что в этот момент вся его агрессия мгновенно испарилась. Вопросы он задавал осторожно и не спеша, ошибки (а их у меня было) поправлял спокойно --- скорее в манере преподавателя, а не экзаменатора. Про последнюю из задач он сказал <<ну, это Вам на дом>>.

Вердикт был примерно такой (как только кончилась математика, к нему снова вернулась агрессивная манера):

--- Странно, что Вы не полный идиот, как кажется вначале. Мы можем Вас взять, но учтите, что если Вы хотите учиться в аспирантуре, не ходить на работу (возможно, он перечислил еще какие-то академические блага), то Вам придется много и тяжело программировать. То, что Вам скажут.

Так уж вышло, что я не выбрал этот вариант (подробности надо вспоминать отдельно). 
\end{quote}
(Конец цитаты.)

\subsection*{ВЗМШ и книжки для неё}

В момент возникновения слово <<ВЗМШ>> было сокращением для <<Всесоюзной Заочной Математической Школы>>, которая была организована по инициативе И.\,М. и поддержке тогдашнего ректора МГУ, замечательного Ивана Георгиевича Петровского в 1964 году (подробнее об истории создания ВЗМШ и о роли различных людей можно прочесть в \url{https://elementy.ru/nauchno-populyarnaya_biblioteka/432911/Zaochnaya_matematicheskaya_shkola}). Попав на первый курс мехмата в 1974, я стал  <<проверяющим>> этой школы и в течении нескольких лет регулярно проверял работы школьников. Школьникам посылались книжки-задания по математике. Точнее, им посылались книжки (изданные типографски или на ротапринте), и прилагался список избранных задач из них, решения которых нужно было записать в тетрадку и прислать на проверку. После проверки (в которой я и участвовал) тетрадки с замечаниями посылались школьникам обратно; если было много ошибок, им предлагалось исправить решения и прислать исправленное для повторной проверки. 

Возвращаясь к книжкам-заданиям: часть из них была к тому времени издана в серии <<Библиотечка физико-математической школы>> в издательстве <<Наука>> (её редактором был Гельфанд), а часть заданий печатались на ротапринте с машинописных оригиналов. (Некоторые из неизданных ротапринтных брошюр я выложил в \url{archive.org}, см.~\url{https://archive.org/details/combinatorics1989}, \url{https://archive.org/details/integer1976}.) Вот изданные в <<Науке>>:
\begin{itemize}
\item И.\,М.\,Гельфанд, Е.\,Г.\,Глаголева, А.\,А.\,Кириллов. Метод координат (4-е издание, 1973, \url{https://math.ru/lib/book/djvu/zaochn/b1_73.djvu})
\item И.\,М.\,Гельфанд, Е.\,Г.\,Глаголева, Э.\,Э.\,Шноль. Функции и графики (3-е издание, 1968, \url{https://math.ru/lib/zaochn/2}).
\item Н.\,Б.\,Васильев, В.\,Л.\,Гутенмахер. Прямые и кривые (2-е издание, 1978, \url{https://math.ru/lib/zaochn/4}).
\item А.\,А.\,Кириллов. Пределы (2-е издание, 1973, \url{https://math.ru/lib/zaochn/7}).
\end{itemize}
Последняя книжка была дополнительным заданием для передовых школьников.
Все они, как я понимаю, были написаны для ВЗМШ; думаю, что и те, в которых И.\,М. не указан как автор, были написаны <<под его присмотром>>. Была ещё книжка <<Уравнения и неравенства>> М.\,И.\,Башмакова (2-е издание, 1976, \url{https://math.ru/lib/book/djvu/zaochn/b5_76.djvu}), которая тоже использовалась как задание и вышла в той же серии, но судя по стилю (и по отсутствию упоминаний в тексте), к ней И.\,М. отношения не имел. 

В серии выходили и другие книжки: <<Задачи по элементарной математике>> (С.\,И.\,Гель\-фанд, М.\,Л.\,Гер\-вер, А.\,А.\,Ки\-рил\-лов, Н.\,Н.\,Кон\-стан\-ти\-нов, А.\,Г.\,Куш\-ни\-рен\-ко, 1965, \url{https://math.ru/lib/book/djvu/zaochn/b3_65.djvu}),  <<Математические задачи>> (Е.\,Б.\,Дын\-кин, С.\,А.\,Мол\-ча\-нов, А.\,Л.\,Ро\-зен\-таль, А.\,К.\,Толпыго, 3-е издание, 1971, \url{https://math.ru/lib/book/djvu/zaochn/d1_71.djvu}), <<Математические соревнования. Арифметика и алгебра>> (Е.\,Б.\,Дын\-кин, С.\,А.\,Мол\-ча\-нов, А.\,Л.\,Ро\-зен\-таль, 1970, \url{https://math.ru/lib/book/djvu/zaochn/d3_70.djvu})
<<Математические соревнования. Геометрия>> (Н.\,Б.\,Ва\-силь\-ев, С.\,А.\,Мол\-ча\-нов, А.\,Л.\,Ро\-зен\-таль, А.\,П.\,Са\-вин,  1974, \url{https://math.ru/lib/book/djvu/zaochn/d4_74.djvu}). Но эти книжки на моей памяти в качестве заданий ВЗМШ не использовались.

Составляя задания для ВЗМШ, авторы столкнулись со сложной задачей. Ученики ВЗМШ --- это школьники, интересующиеся математикой (иначе бы они не стали  писать вступительную работу в ВЗМШ). Соответственно, с школьными уроками и учебниками математики они были в той или иной степени знакомы, и повторять их, давая типовые, пусть даже и более трудные задачи, не имело смысла. С другой стороны, даже и интересующиеся математикой школьники в среднем понимали в школьном курсе математики не так много, так что опираться на изученное в школе, считая его известным, тоже было нельзя. Можно было бы пытаться выйти из положения, давая задачи на темы, мало представленные в школьной программе (<<олимпиадные задачи>>, <<задачи на сообразительность>>) ---  но всё-таки для \emph{школы}, пусть и заочной, этого мало, невозможно не пересечься со школьной программой. Значит, надо было рассказывать о том же, но не так же, найти какую-то свежую точку зрения, новые неожиданные вопросы (которые в процессе обучения зачастую важнее, чем ответы), <<увидеть всё заново>>, если пользоваться выражением Честертона.

Преподаватели математики, вынужденные следовать программе (будь то школьное или высшее образование), обнаруживают себя в роли машиниста паровоза, которому велено следовать расписанию --- и не обращать внимания на то, что большинство вагонов уже давно отцепились. При этом надо поддерживать впечатление (иллюзию) учебного процесса --- когда происходят лекции, семинары, экзамены, ставятся оценки и т.п. Результатом часто становится, как говорят, <<плохое равновесие>> --- в экзамены включаются задания, которые можно выполнять по образцу, не понимая существа дела,\footnote{Байки от И.\,М.:
\begin{itemize}
\item Занятие в вечерней школе (<<рабочей молодёжи>>): <<Что больше, $3/2$ или $5/3$?>> --- Молчание. Отчаявшийся преподаватель: <<Ну сами подумайте: что лучше, три бутылки на двоих, или пять на троих?>> --- Слушатели, недоумевая: <<Конечно, пять на троих лучше, но при чём тут это?>>
 
\item Студент мехмата, приехавший на каникулы и зашедший в родную школу, на перемене учительнице: <<Вы знаете, Марья Ивановна, складывать дроби отдельно по числителям и знаменателям неправильно, надо\ldots>> (объясняет) --- Марья Ивановна, после перемены: <<Дети, из Москвы пришло новое указание\ldots>>
\end{itemize}
Ещё одну историю (уже не от Гельфанда) я прочёл давно в какой-то книжке, но не могу сейчас найти источник и пересказываю по памяти: рассказчик, ставший потом математиком, вспоминает, что в гимназии его товарищ плохо отвечал на вопросы, и в конце концов учитель (физики?) решить спросить что-нибудь простое, и когда ученик на вопрос, куда покатится шар, положенный на наклонную плоскость, ответил, что вверх --- выгнал его со скандалом, решив, что ученик над ним издевается. Рассказчик потом спросил своего товарища: <<ты что, действительно не знал, что шар покатится вниз?>> --- <<\emph{Настоящий} шар, конечно, покатится вниз, но\ldots>>

Плохое равновесие можно было наблюдать и на советских вступительных экзаменах и <<задачах из Сканави>> с бесконечными ОДЗ --- а теперь оно с книжками типа <<как решать задачу B7 на ЕГЭ>>. } 
 а занятия превращаются в демонстрацию образцов. Вернуться же за вагонами --- прямой путь к авариям,\footnote{Кажется, примерно это случилось с Тоомом, когда он пытался преподавать в американских университетах <<второй лиги>> (а скорее третьей) --- до того, как переехал в Бразилию. Он пытался реально кого-то чему-то научить, но выяснилось, что никому это не надо, а только раздражает.} да и сами вагоны уже не очень годятся для движения.

При такой общей обстановке задача ВЗМШ кажется почти безнадёжной: задания одновременно должны были быть и посильными для большинства, и интересными для самых продвинутых школьников, одновременно простыми и нетривиальными, школьными и не-школьными. Удивительно, в какой степени авторам удалось выполнить эти противоречивые требования. В полной мере это можно оценить, лишь работая со школьниками (в том числе по этим книжкам), но всё же расскажу о книжках подробнее и приведу несколько примеров, чтобы было понятно, о чём речь (а также несколько примеров того, как --- на мой взгляд --- делать не надо).

\clearpage
\subsection*{О  конкретных книжках ВЗМШ}

\subsubsection*{Метод координат} 

\begin{wrapfigure}{r}{0pt}
\includegraphics[width=0.3\textwidth]{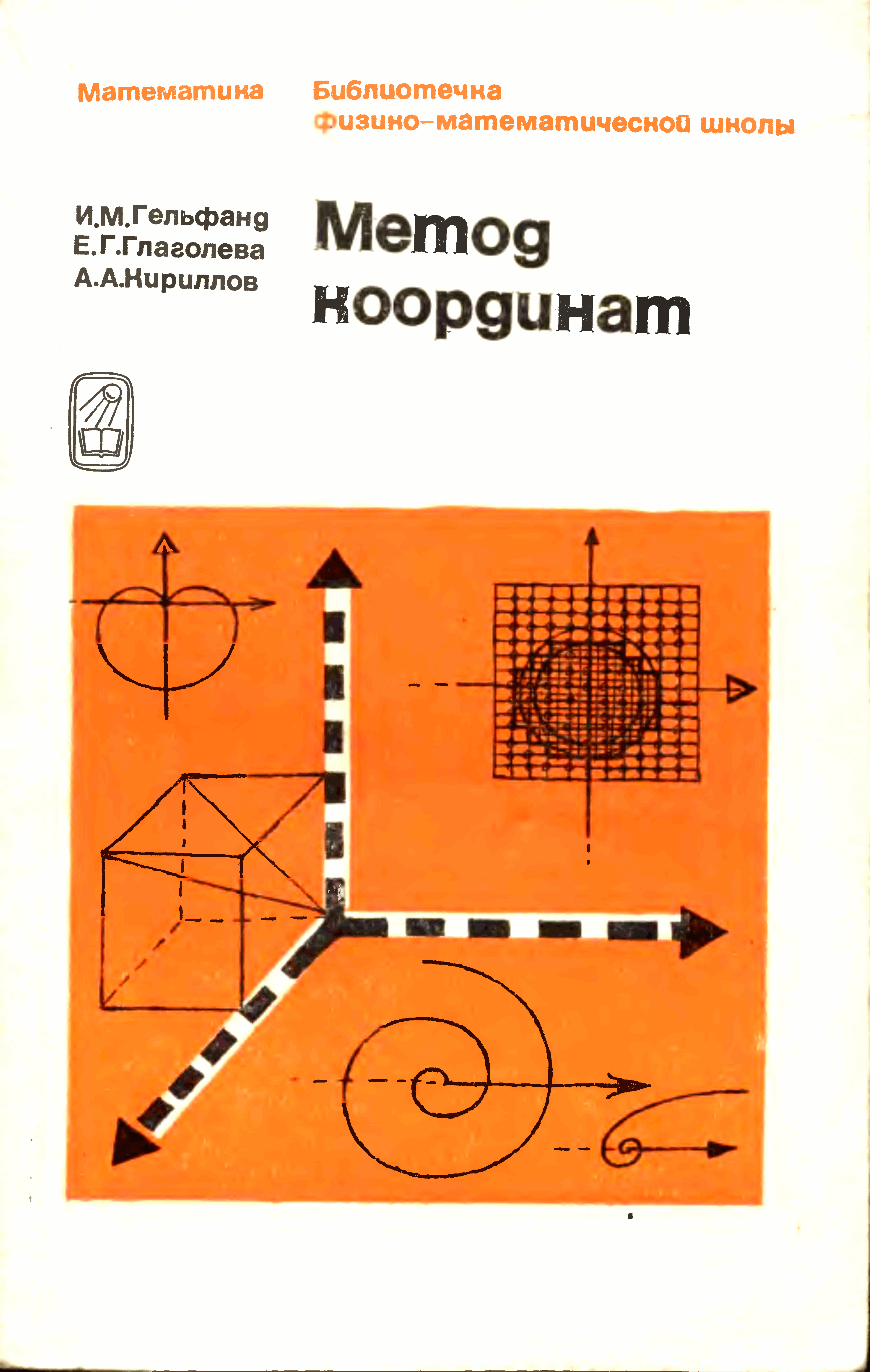}
\end{wrapfigure}

Книжка состоит из разделов, идущих в порядке возрастания размерности: координаты на прямой, на плоскости, в пространстве, и, наконец в четырёхмерном пространстве (вторая половина книжки). 

Что, казалось бы, можно такого неожиданного спросить про <<числовую прямую>> после школьных учебников? Оказывается, что совсем простые вопросы, если они хотя бы немного непривычные, у многих школьников вызывают сложности.  Чему равно $|{-a}|$, если $a$ --- отрицательное число: $a$ или $-a$?  Как (знакомое многим школьниками того времени) уравнение типа $|x+1|+|x+2|=2$ можно решить в уме, если понимать, что абсолютная величина означает расстояние?

Раздел про координаты на плоскости начинается с совсем простых упражнений: не рисуя точки $A(1,-3)$ на плоскости, скажите, в какой четверти она расположена? расшифруйте слово, заданное координатами точек, образующих буквы этого слова. Как найти координату четвёртой вершины параллелограмма с вершинами в $(0,0)$, $(x_1,y_1)$, $(x_2,y_2)$?  Где расположены точки, у которых одна из координат (абсцисса или ордината) равна $3$? А в конце раздела обсуждается, почему множество точек $M$ с данным отношением расстояний $|MA|/|MB|$ до двух фиксированных точек $A$ и $B$ является окружностью и как простая картинка позволяет в уме определить, при каких значениях $a$ система уравнений
    $$
    \left\{
    \begin{aligned}
        x+y=a\\
        x^2+y^2=1
    \end{aligned}\right.
    $$
имеет ровно одно решение.

Координаты в пространстве мало используются в школе (и мало использовались даже в то время, когда были в ходу задачи по стереометрии), так что эта часть мало пересекается с школьной программой. Подробно обсуждаются простые примеры, которые обычный учебник бы скорее всего даже не упомянул, скажем, почему условие $x^2+y^2=1$ задаёт не окружность, а цилиндр, или каким условиям удовлетворяют координаты точек, лежащих внутри единичного куба.

Переходя ко второй части, где речь идёт о четырёхмерном пространстве, авторы мотивируют этот переход таким примером: количество решений неравенства $x^2+y^2\le n$ в целых числах можно примерно найти, зная площадь круга радиуса $\sqrt{n}$, аналогично для $x^2+y^2+z^2\le n$, но что делать для $x^2+y^2+z^2+t^2\le n$? Дальше они предлагают читателю самому попытаться понять, как устроен четырёхмерный куб. Это редкая (и, пожалуй, почти что уникальная для школьников) ситуация, где надо не что-то вычислить или что-то доказать, а придумать естественное определение вместе с авторами книжки, и научиться представлять себе определяемые таким образом объекты. Вот характерная фраза, отражающая стиль всей книжки: 

\begin{quote}
Не надо огорчаться, что мы не привели пока рисунок четырёхмерного куба --- мы это сделаем потом (не удивляйтесь, что можно нарисовать четырёхмерный куб: ведь рисуем же мы трёхмерный куб на плоском листе бумаги). Для этого сначала надо разобраться, как этот куб <<устроен>>, какие элементы в нём можно различать. 
\end{quote}

\subsubsection*{Функции и графики}

\begin{wrapfigure}[13]{r}{0pt}
\includegraphics[width=0.3\textwidth]{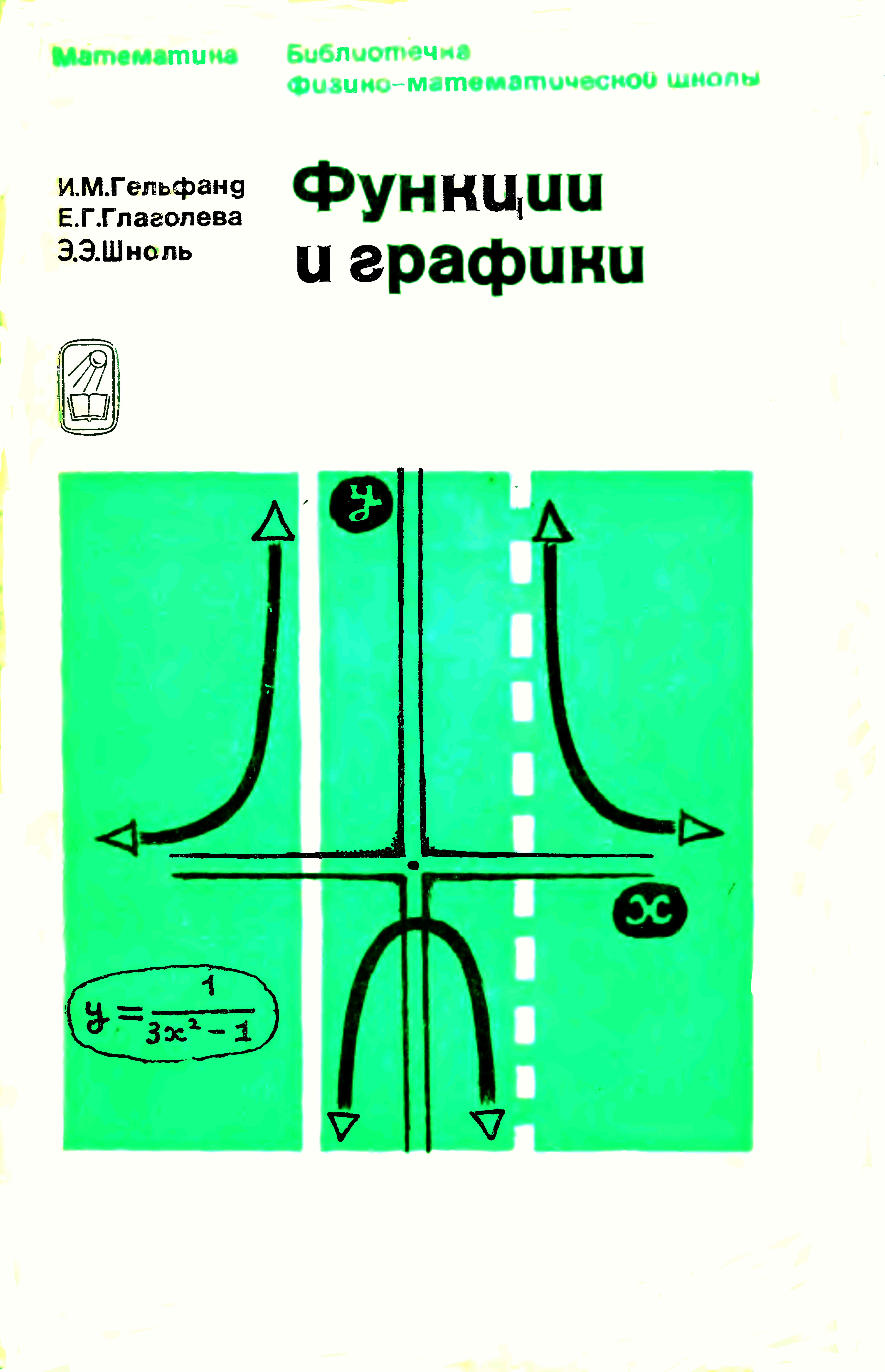}
\end{wrapfigure}

Графики функций --- тоже стандартная для школы тема: сначала рисуют графики линейных функций, потом параболы (графики квадратных трёхчленов). И тем не менее для большинства школьников остаётся довольно загадочным, что это за функции, какие у них бывают графики, и зачем всё это надо. А уж если начать объяснять, что функция --- это множество пар с определёнными свойствами\ldots

Авторы начинают с другого: рисуют кардиограмму, данные сейсмографа и характеристику диода, и не боятся использовать оборот \emph{это когда}, от которого в школе старательно отучают: <<функция --- это когда каждому значению некоторой величины, которую математики называют аргументом и обозначают обычно буквой $x$, отвечает значение другой величины $y$, называемой функцией. $\langle\ldots\rangle$ Сила тока в полупроводниковом элементе есть функция напряжения, так как каждому значению напряжения соответствует определённое значение силы тока.>>

Дальнейшее изложение также совсем не напоминает школьную дрессировку (как рисовать график линейной функции, затем квадратной и пр.): всё это тоже будет, но одним из первых примеров служит график функции $x\mapsto \lfloor x\rfloor$ (целая часть), а обучение построению графиков начинается с разбора графиков $y=1/(1+x^2)$ и $y=1/(3x^2-1)$, которые строятся по точкам с подробным обсуждением, какие новые особенности становятся видны при добавлении точек. 

Дальше идёт более стандартный материал: про график линейной функции, график модуля. И в качестве, как теперь говорят, <<бонуса>> предлагается задача, в которой условие не содержит никаких упоминаний функций и графиков, но естественное решение их использует:
\begin{quote}
На окружности расположено $7$ коробок со спичками. В первой лежит $19$ спичек, во второй $9$, в остальных соответственно $26,8,18,11,14$ спичек. Разрешается перекладывать спички из любой коробки в любую из соседних с ней. Требуется переложить спички так, чтобы во всех коробках их стало поровну. Как это сделать, перекладывая как можно меньше спичек?
\end{quote}

Вообще один из важных принципов, последовательно проводимых в этой и других книжках --- в каждый момент рассказа должно быть понятно, к чему мы стремимся, что сейчас делаем, какой план дальнейших действий и пр. Эта структура, хотя и никак не отражается в формальной модели математики как логической теории (где любое допустимое по правилам вывода действие одинаково допустимо), является ключевой для читателя. К сожалению, во многих математических статьях она почти что отсутствует. Но особенно плохо это в книжках для школьников: увы, часто их авторы отвечает на вопрос, который у школьника даже и не возник. 

Гельфанд с соавторами этого избегают.  Скажем, стандартный приём выделения полного квадрата появляется здесь не просто так, а с заранее понятной целью: показать, что график $y=x^2+px+q$ <<по форме ничем не отличается от параболы $y=x^2$ и лишь занимает другое положение относительно координатных осей>>. (Обратите внимание, кстати, на бытовое <<по форме ничем не отличается>> вместо педантичного <<получается параллельным переносом>>.)

А (в следующем разделе) построение графика гиперболы $y=1/x$ начинается так: 

\begin{center}
\includegraphics[width=0.7\textwidth]{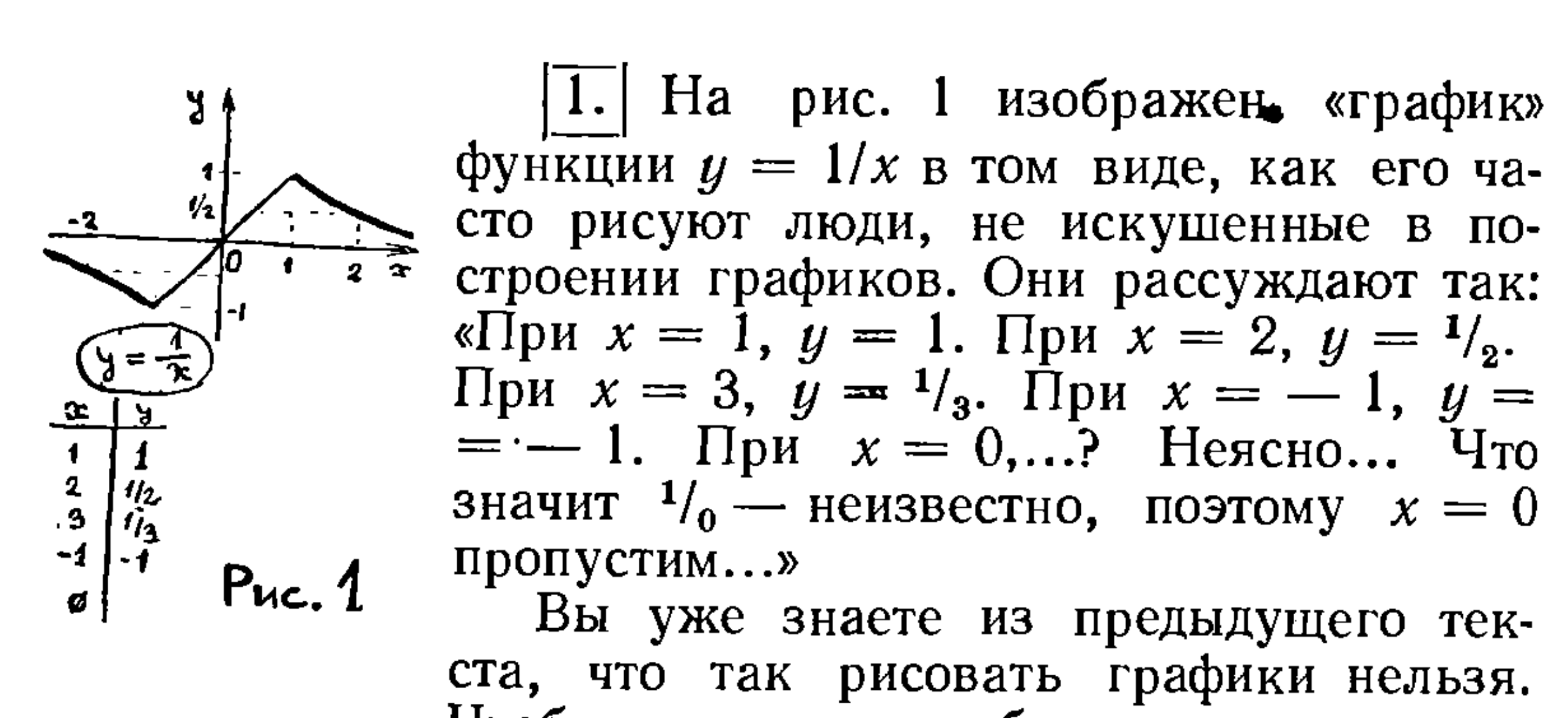}
\end{center}

Ещё один принцип, которому следуют авторы: о чём бы они не рассказывали, это всегда остаётся частью большой картины, и по ходу дела проявляются связи с другими её частями. Не отступая далеко от основного изложения, они замечают, что построение графиков функций типа $x\mapsto f(x+a)$ или $x\mapsto -f(x)$ напоминает о движениях плоскости, линейная функция связана с равномерным движением и арифметическими прогрессиями, прямая может проходить только через одну целочисленную точку, потому что бывают иррациональные числа, разности последовательных значений квадратного трёхчлена в целых точках образуют арифметическую прогрессию, парабола --- не только график квадратного трёхчлена, но и множество точек, равноудалённых от фокуса и директрисы, а также сечение конуса, траектория падающего тела и форма воды во вращающемся стакане, кратный корень квадратного уравнения превращается в пару корней при малом шевелении параметра, графики функций имеют асимптоты,  разные функции растут с разной скоростью, касательную лучше определять как предел секущей, а не как прямую, имеющую одну общую точку, и многое другое.

\clearpage
\subsubsection*{Прямые и кривые}
\begin{wrapfigure}[13]{r}{0pt}
\includegraphics[width=0.301\textwidth]{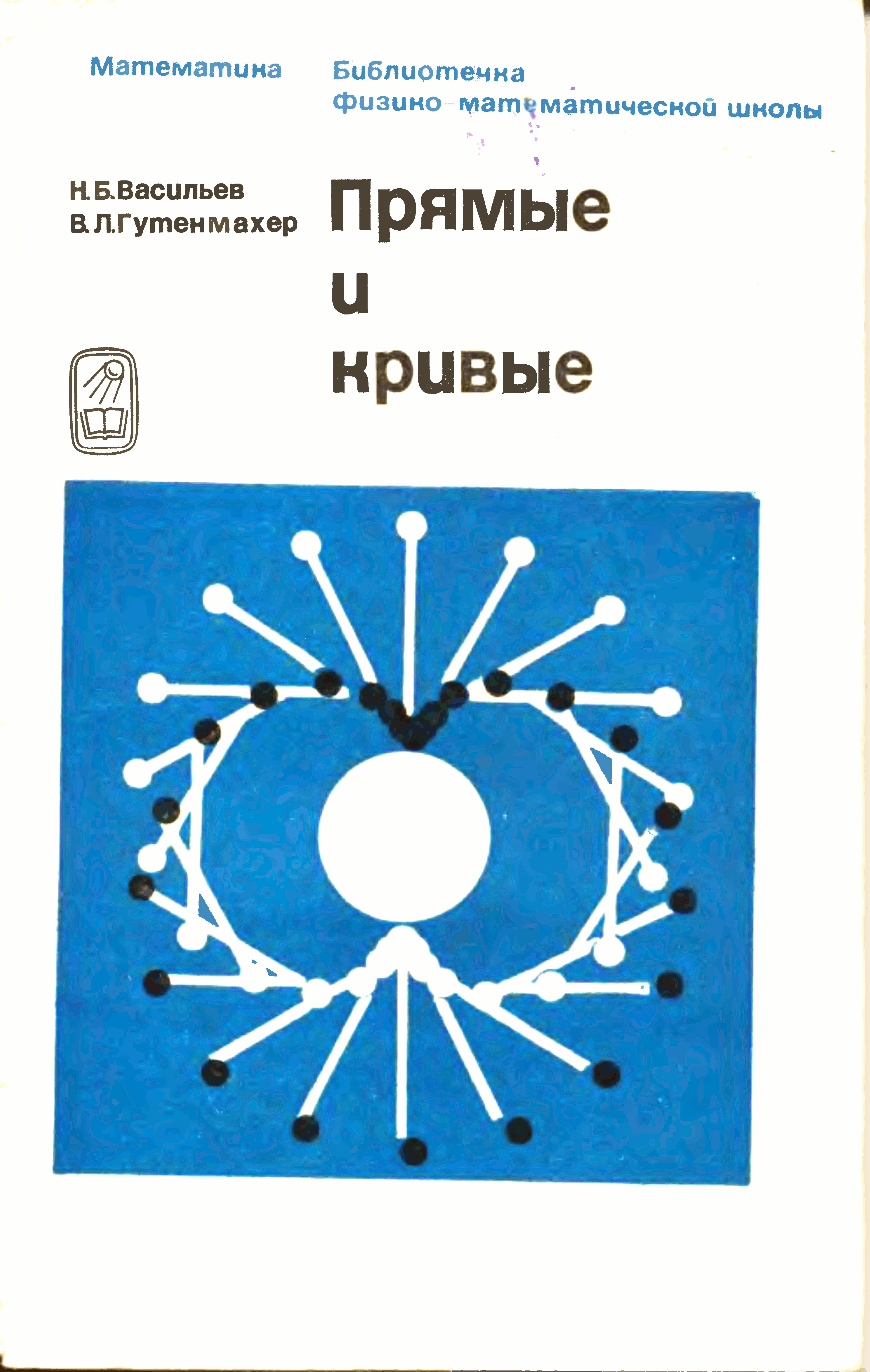}
\end{wrapfigure}

Геометрическая часть программы ВЗМШ была делом особенно сложным: здесь тоже надо было найти какой-то новый взгляд на предмет, при наличии двухтысячелетней традиции. Тем не менее авторам геометрической книжки (Н.\,Б.\,Ва\-силь\-ев и В.\,Л.\,Гу\-тен\-ма\-хер --- для их совместных текстов даже был псевдоним <<Вагутен>>) это в значительной мере удалось. Книжка эта представляет собой <<тему с вариациями>>, и темой являются <<геометрические места точек>> ---  множества точек, задаваемые некоторыми условиями. 

\bigskip
Сами авторы в первом же абзаце описывают идею книжки так: 
\begin{quote}
Главные действующие лица этой книжки --- различные геометрические фигуры, или, как они здесь чаще называются, <<множества точек>>. Вначале появляются самые простые фигуры в различных сочетаниях. Они двигаются, обнаруживают новые свойства, пересекаются, объединяются, образуют целые семейства и меняют своё обличье --- иногда до неузнаваемости; впрочем, интересно увидеть старых знакомых в сложной обстановке, в окружении новых фигур, появляющихся в финале.
\end{quote}

Вот несколько примеров, показывающих, что имеется в виду.

Кошка сидит на середине лестницы $KL$, которая скользит у стены; по какой траектории она движется? (Медиана прямоугольного треугольника, проведённая к гипотенузе, равна её половине.)

\begin{center}
\includegraphics[width=0.3\textwidth]{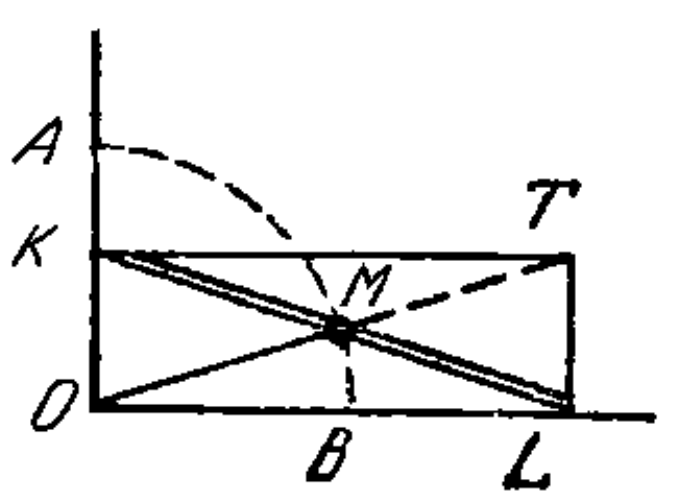}
\end{center}

А если взять другую точку, то получится эллипс, и так устроен эллипсограф Леонардо да Винчи.

Другая задача о движении: какую траекторию описывает точка окружности, катящейся без проскальзывания внутри вдвое большей неподвижной окружности? Оказывается, что по прямой, и это называется <<теорема Коперника>> и следует из теоремы о вписанном угле школьного курса геометрии.

\begin{center}
\includegraphics[width=0.3\textwidth]{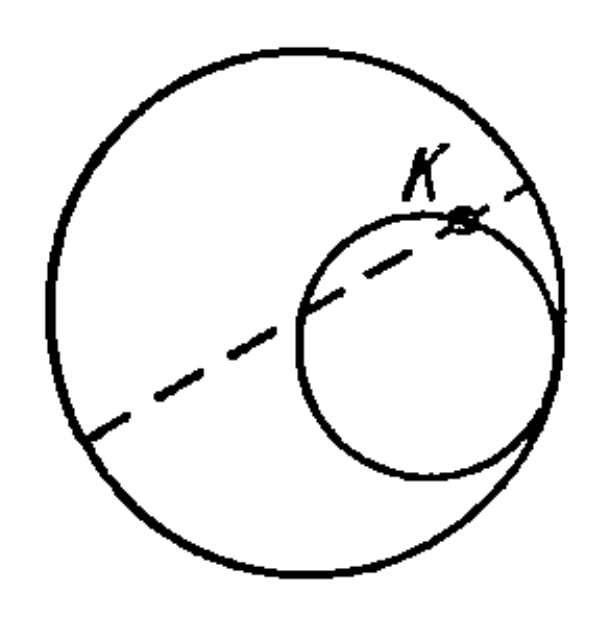}
\end{center}

Сама эта теорема тоже может быть сформулирована наглядно:

\begin{center}
\includegraphics[width=0.7\textwidth]{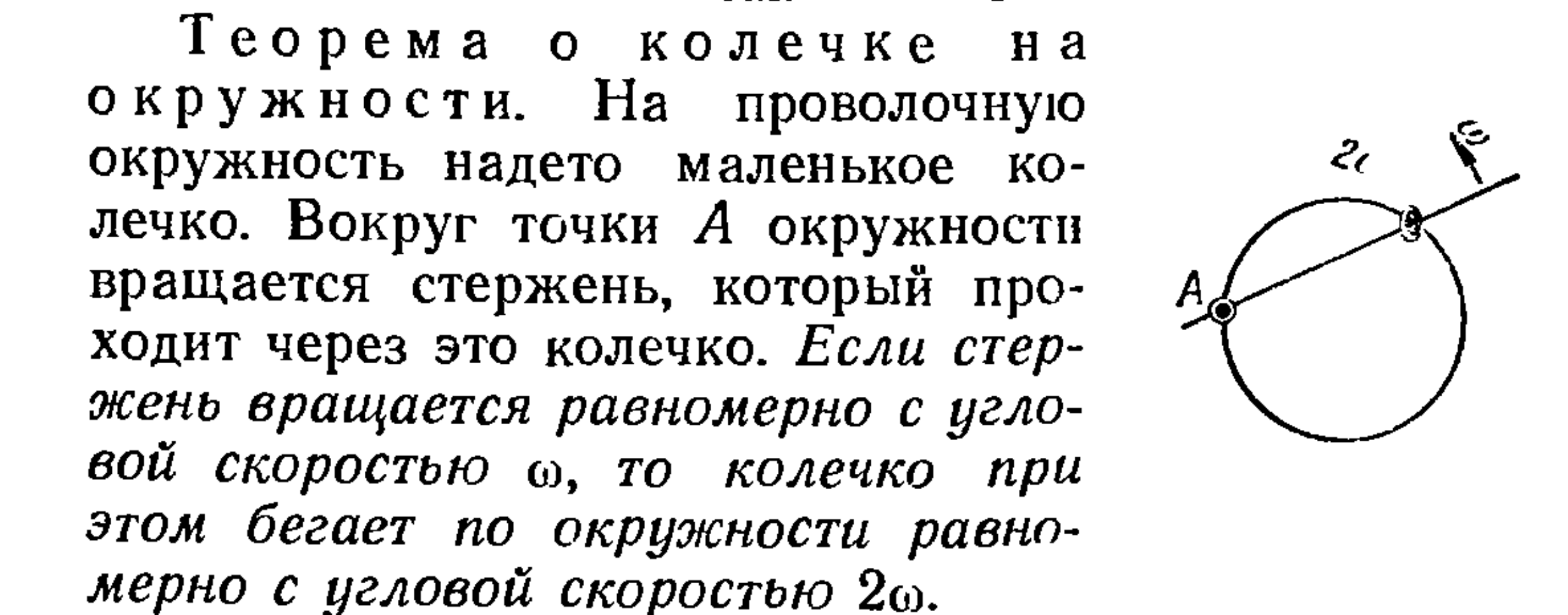}
\end{center}

Другая кинематическая интерпретация: 

\begin{center}
\includegraphics[width=0.7\textwidth]{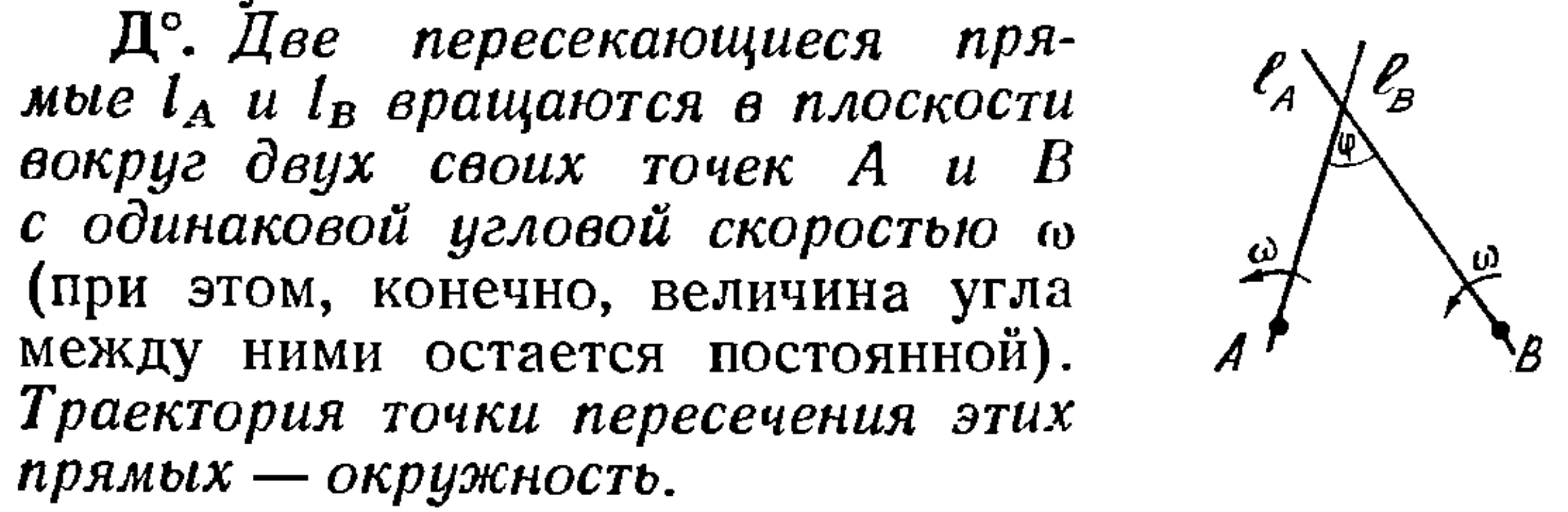}
\end{center}

Наконец, в конце книжки та же самая теорема появляется в ещё одной ситуации, когда решается такая задача:

\begin{center}
\includegraphics[width=0.7\textwidth]{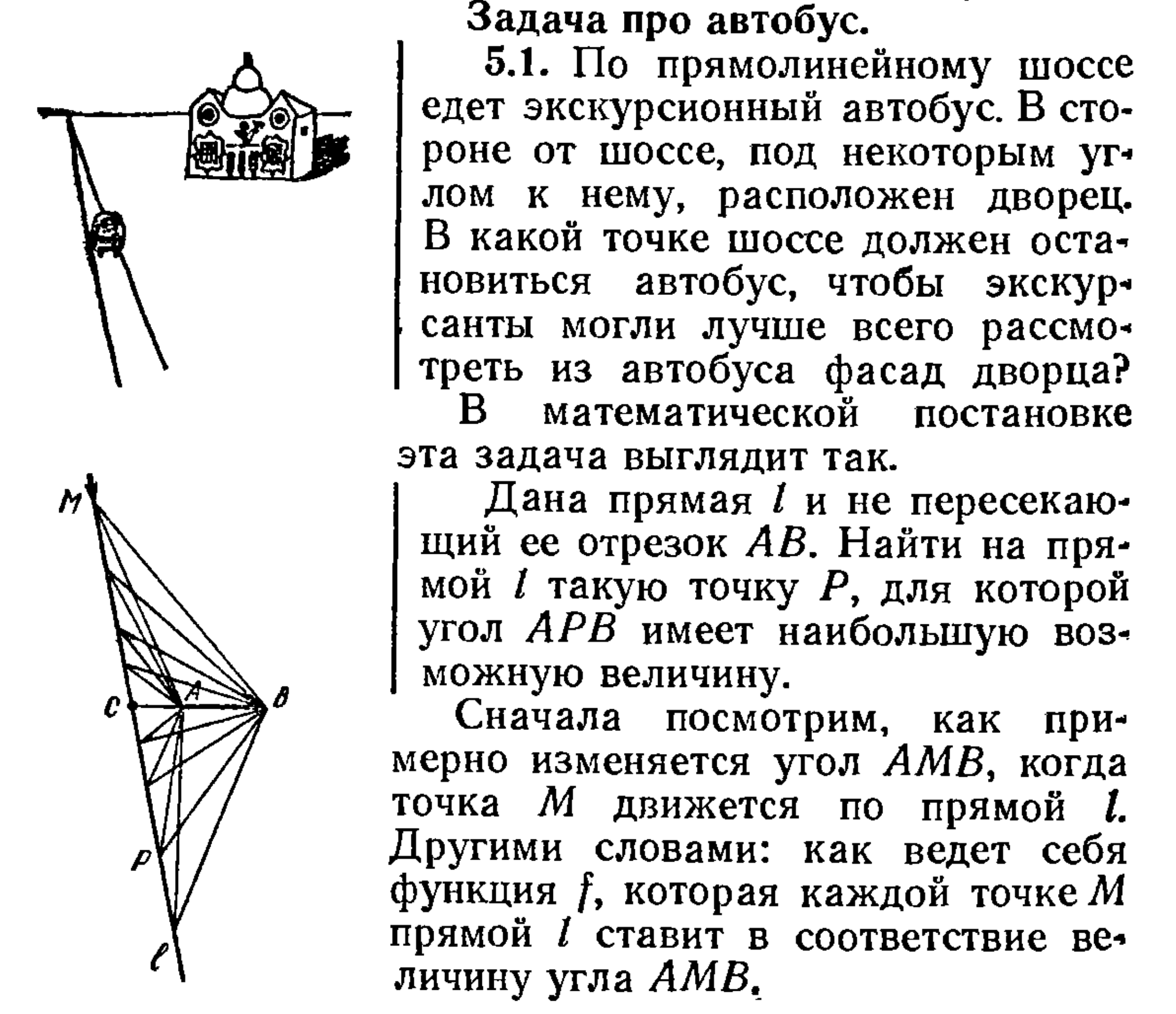}
\end{center}

Чтобы решить её, посмотрим на поведение угла как функции точки $P$, освободив точку $P$ от ограничения <<лежать на прямой $l$>>. Это будет функция на плоскости, и по теореме о вписанном угле её линиями уровня будут дуги окружностей. Нам нужно найти линию максимального уровня, имеющую общую точку с прямой $l$, она будет окружностью, касающейся прямой $l$ и проходящей через точки $A$ и $B$:

\begin{center}
\includegraphics[width=0.3\textwidth]{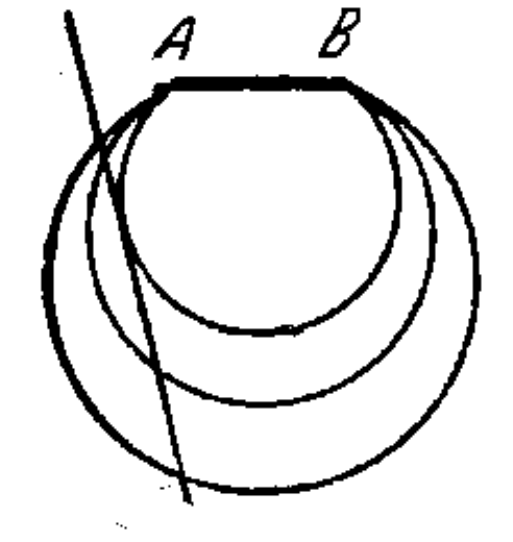}
\end{center}

Точнее надо говорить о дуге такой окружности, и таких дуг две (точка касания может быть выше или ниже $AB$). Соответствующие углы надо сравнить, и выбрать ту точку, где угол больше (или ту, откуда виден фасад дворца, если вернуться к исходной постановке). Таким образом, задача об автобусе сводится к классической задаче на построение (провести окружность, касающуюся данной прямой, через две данные точки).

\clearpage
\subsubsection*{Задачи заочных математических олимпиад}

\begin{wrapfigure}[9]{r}{0pt}
\includegraphics[width=0.35\textwidth]{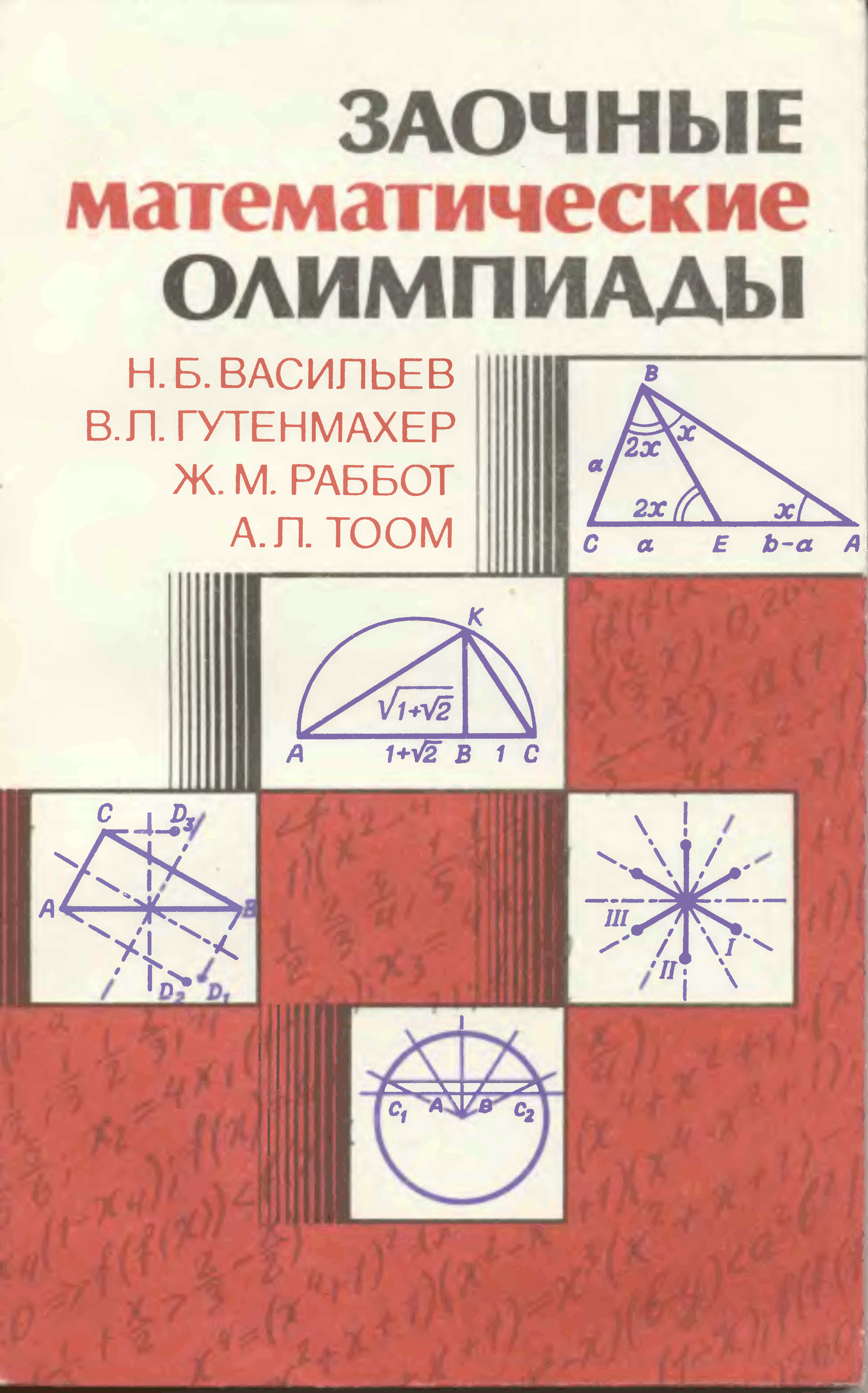}
\end{wrapfigure}

Ещё одна традиция, сложившаяся в заочной школе, связана с задачами для вступительных работ. В лучшие  годы среди поступающих был конкурс, и для поступления нужно было пройти отбор, написав вступительную работу. Работа эта была заочной (в моё время условия задач печатались в журнале <<Квант>>, и школьники могли присылать решения по указанному там адресу). 
\vspace*{4.5cm}
\medskip

Подбор задач для такой <<заочной олимпиады>> --- дело непростое, они должны удовлетворять сразу нескольким требованиям:
\begin{itemize}
\item задача должна быть не типовой, иначе школьникам будет неинтересно;
\item задача должна быть несложной, иначе мало кто её решит, ведь решающие её только поступают в школу;
\item задача (если её решить правильно) должна быть несложной для записи, потому что никакого опыта записи решений задач, не укладывающихся в школьный стандарт, у большинства школьников нет;
\item задача должна быть простой для проверки, потому что работ много и в проверке участвуют преподаватели (студенты мехмата) без большого опыта чтения работ школьников;
\item наконец, желательно минимизировать вероятность ситуации, когда школьнику кажется, что задачу он решил правильно, а её не зачли (в частности, задачи <<на доказательство>> подходят плохо).
\end{itemize}

Постепенно в ВЗМШ собралась коллекция подходящих задач, и сначала сборник таких задач был издан как брошюра общества <<Знание>>, а потом вышла и книжка (<<Заочные математические олимпиады>>, авторы Н.\,Б.\,Васильев, В.\,Л.\,Гутенмахер, Ж.\,М.\,Раббот, А.\,Л.\,Тоом, издательство <<Наука>>, 1987, второе издание).
  
Вот несколько примеров (как мне кажется) удачных задач из этой книжки:

\begin{itemize}
\item Какое наибольшее количество воскресений может быть в году?
\item Существует ли такое целое число, которое при зачёркивании первой цифры уменьшается в $57$~раз?
\item Точки $A$, $B$, $C$ являются вершинами неравнобедренного треугольника. Сколькими способами можно поставить на плоскости точку $D$, чтобы множество точек $\{A,B,C,D\}$ имело ось симметрии?
\item Какое наименьшее число участников может быть в математическом кружке, если известно, что девочки составляют в нём меньше $50\%$, но больше $40\%$?
\item Известно, что доля блондинов среди голубоглазых больше, чем доля блондинов среди всех людей. Что больше: доля голубоглазых среди блондинов или доля голубоглазых среди всех людей?
\item Поезд двигался в одном направлении $5{,}5$ часов. Известно, что за любой отрезок времени длительностью в один час он проезжал ровно $100$~км. Верно ли, что поезд ехал равномерно? Верно ли, что средняя скорость поезда равна $100$~км/ч?
\end{itemize}  

\clearpage
\subsubsection*{Тригонометрия}

Как я уже упоминал, И.\,М. написал две книжки по тригонометрии с разными соавторами: есть русская книжка (c С.\,М.\,Львовским и А.\,Л.\,Тоомом, МЦНМО) и американская (соавтор Mark Saul, Birkh\"auser). 

\begin{center}

\includegraphics[width=0.3\textwidth]{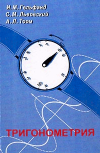}\qquad	\qquad
\includegraphics[width=0.3\textwidth]{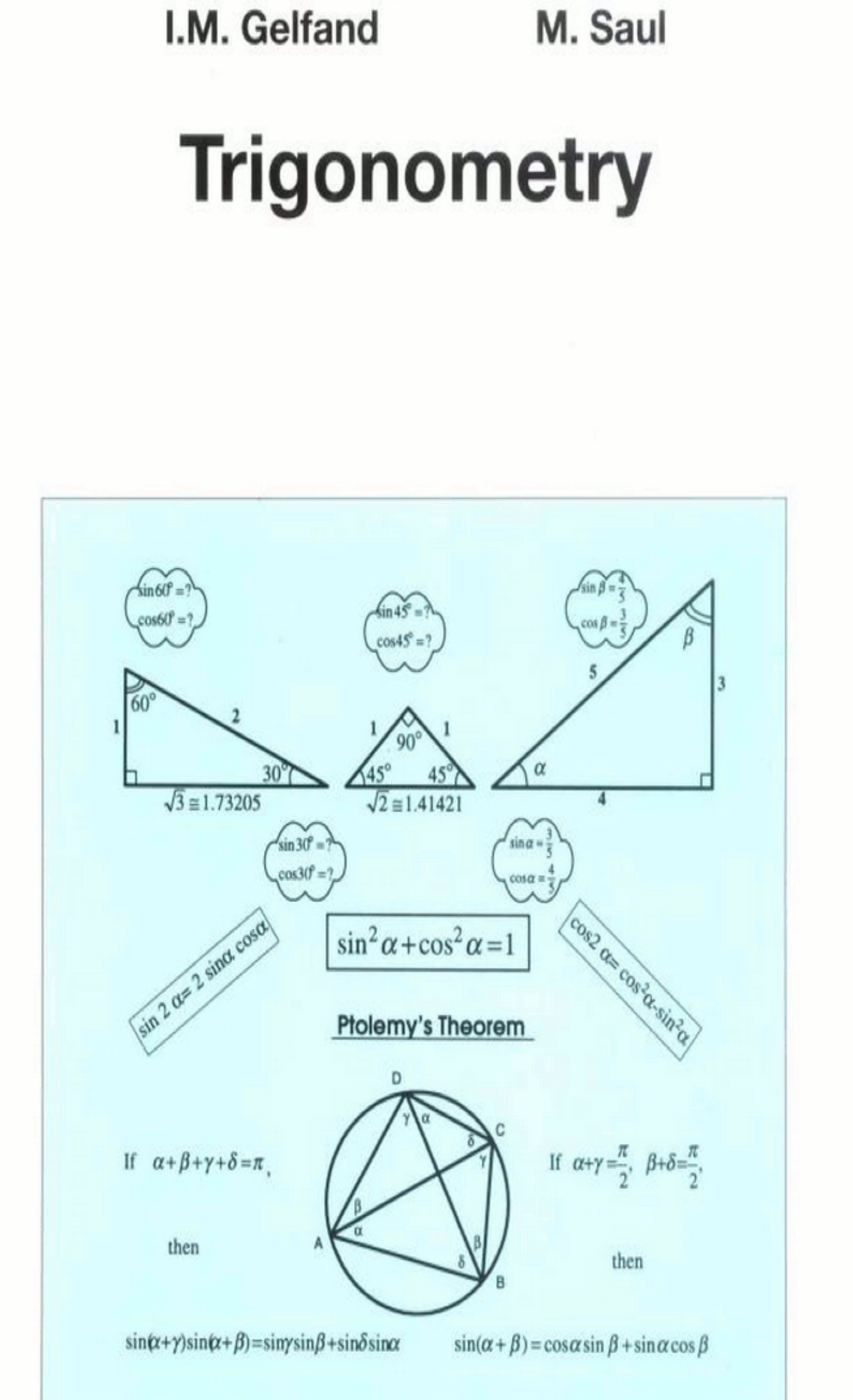}\\
(Рисунок на обложке американской книжки сделан Т.\,Алексеевской)
\end{center}

Как-то в разговоре И.\,М. сказал, что одну и ту же книжку с разными людьми он писал бы по-разному. И действительно, эти книжки отличаются не только по языку, но и по стилю и содержанию (отчасти и потому, что аудитории различаются: американским школьникам, например, ни к чему приёмы решения тригонометрических уравнений вступительного типа).

Но и в американских книжках основные принципы остаются теми же. Вот несколько примеров. Так объясняется, в каких случаях график сдвигается по горизонтали:
\begin{center}
\includegraphics[width=0.7\textwidth]{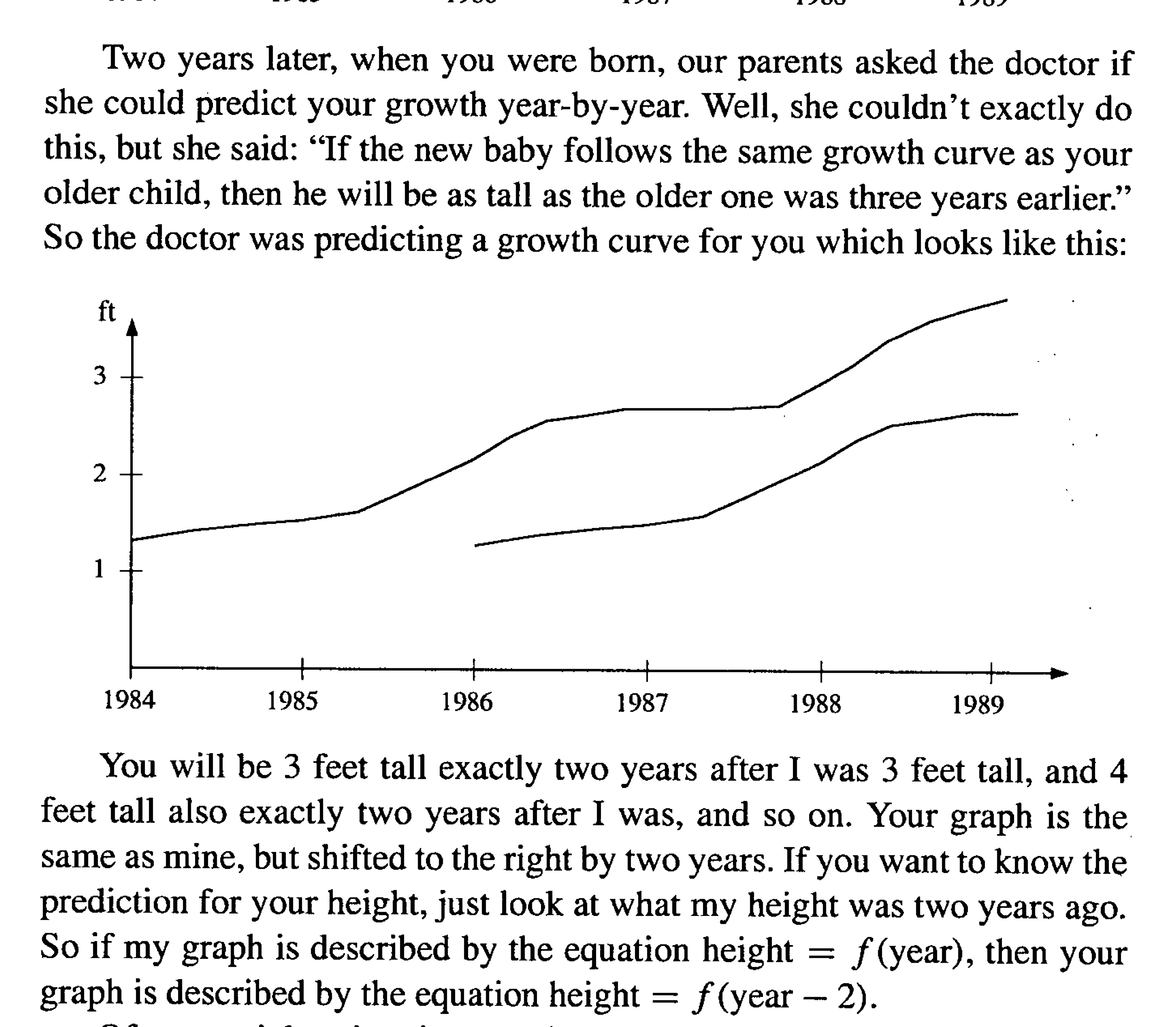}
\end{center}

А вот объяснение, почему формула Герона и должна иметь такой странный на первый взгляд вид:
\begin{center}
\includegraphics[width=0.7\textwidth]{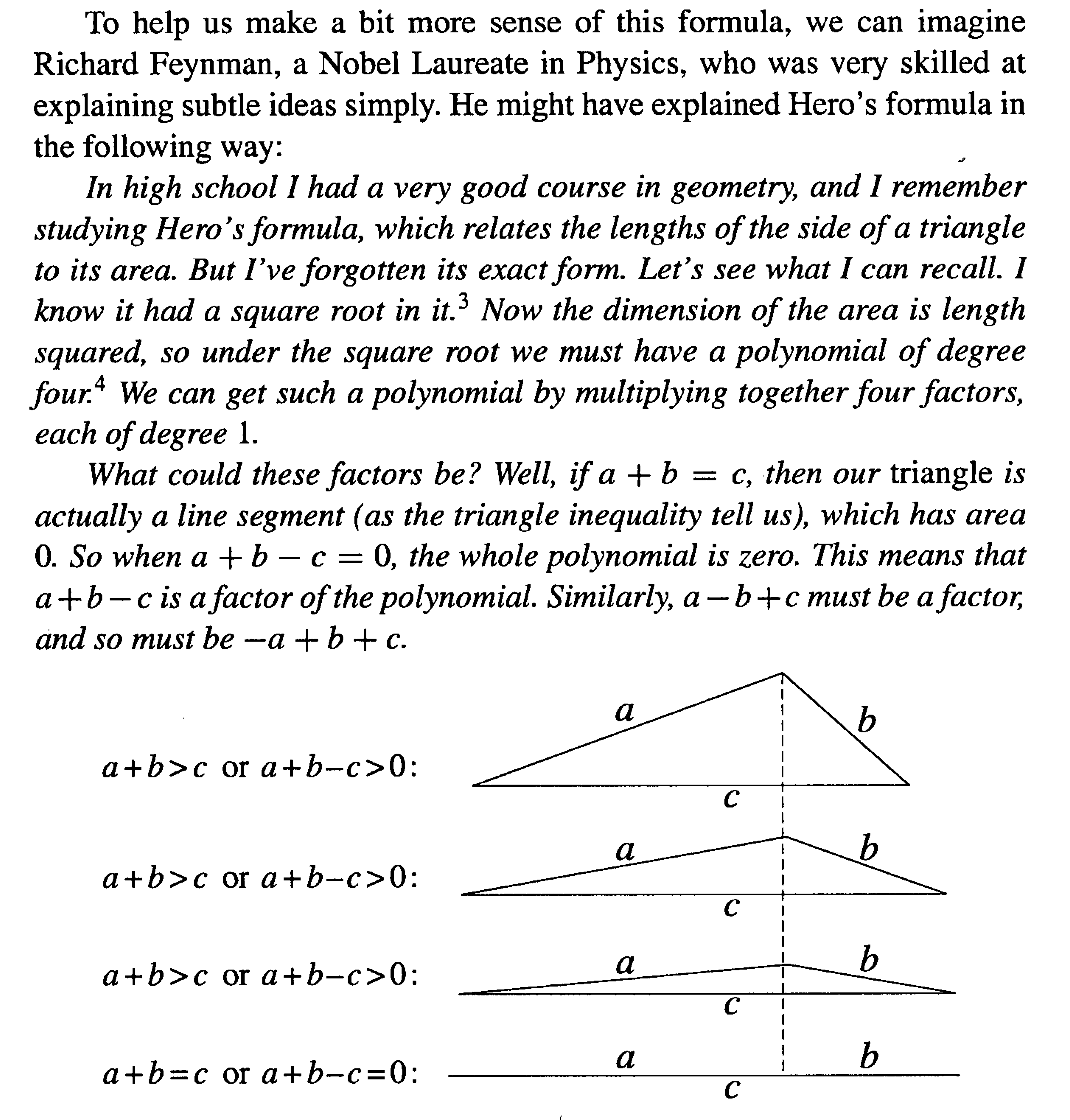}\\
\includegraphics[width=0.7\textwidth]{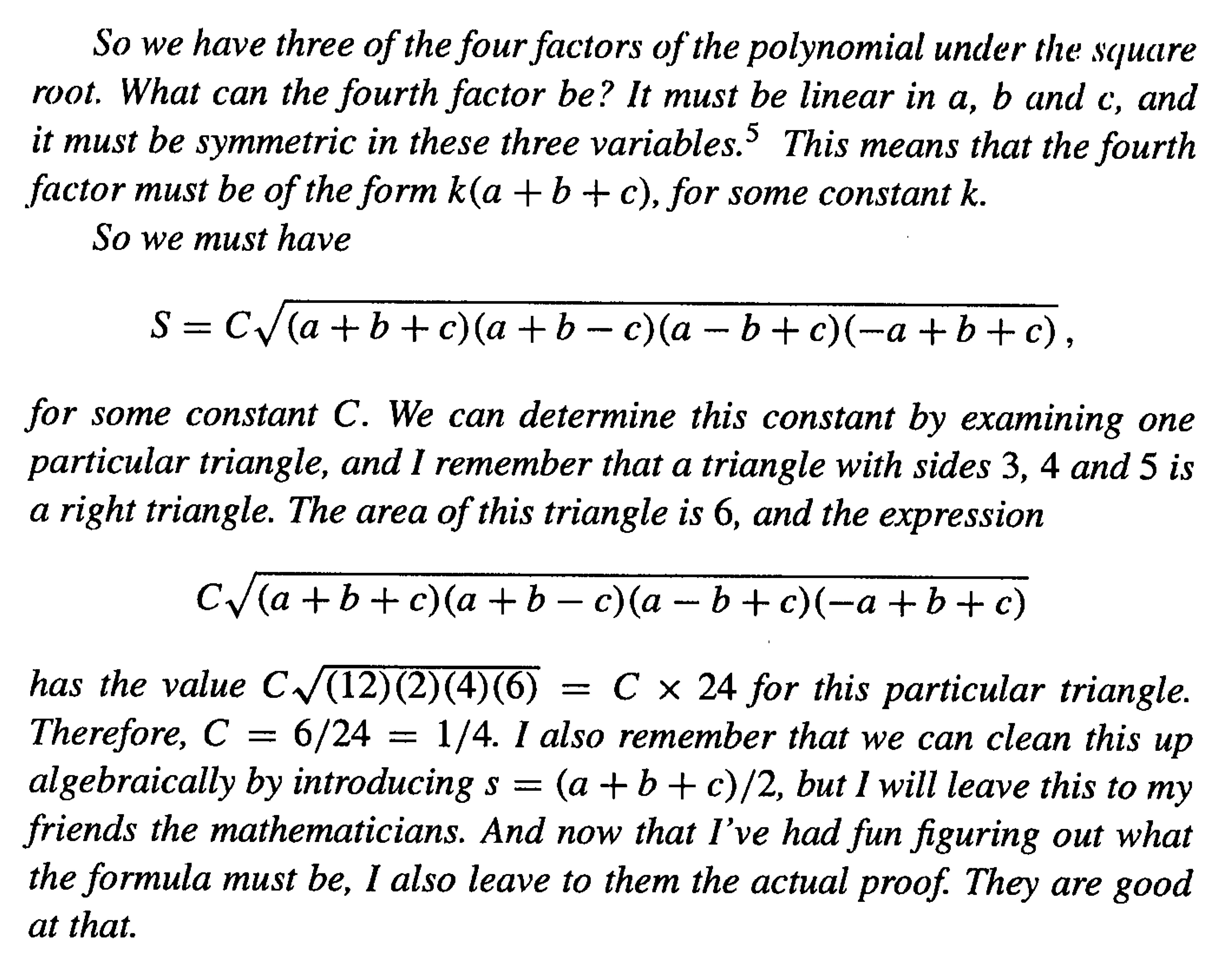}
\end{center}
То, что это объяснение как бы приписывается Фейнману, пожалуй, <<на грани фола>> (и наверняка это Гельфанд, а не Сол), но в данном конкретном случае, думаю, Фейнман был бы только рад.

Пример графика --- не только иллюстрация, но и подготовка к разложению меандра в ряд Фурье:
\begin{center}
\includegraphics[width=0.5\textwidth]{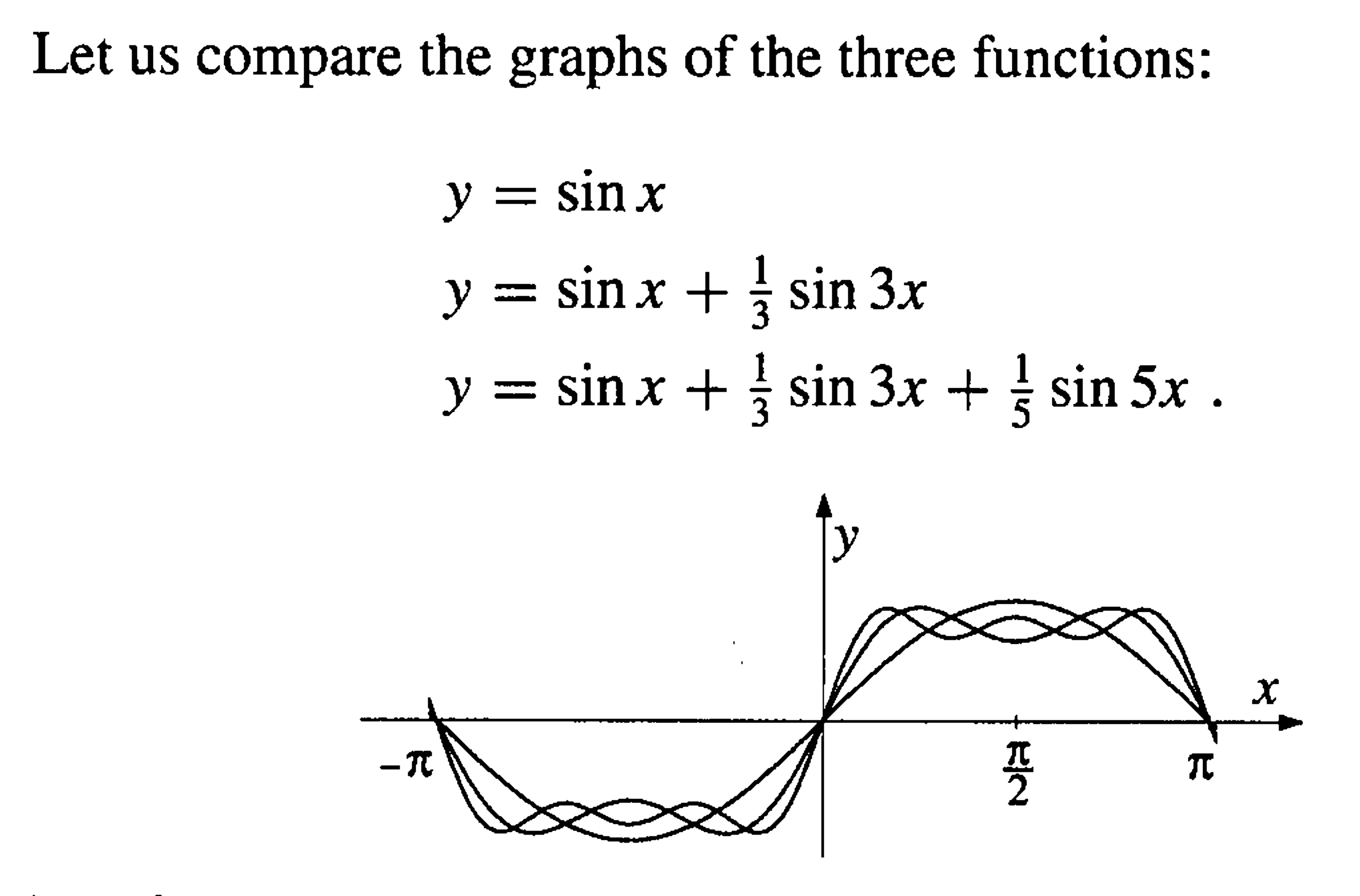}
\end{center}

А вот пример красивой задачи по геометрии, одновременно иллюстрирующей вычисления с тангенсами:

\begin{center}
\includegraphics[width=0.7\textwidth]{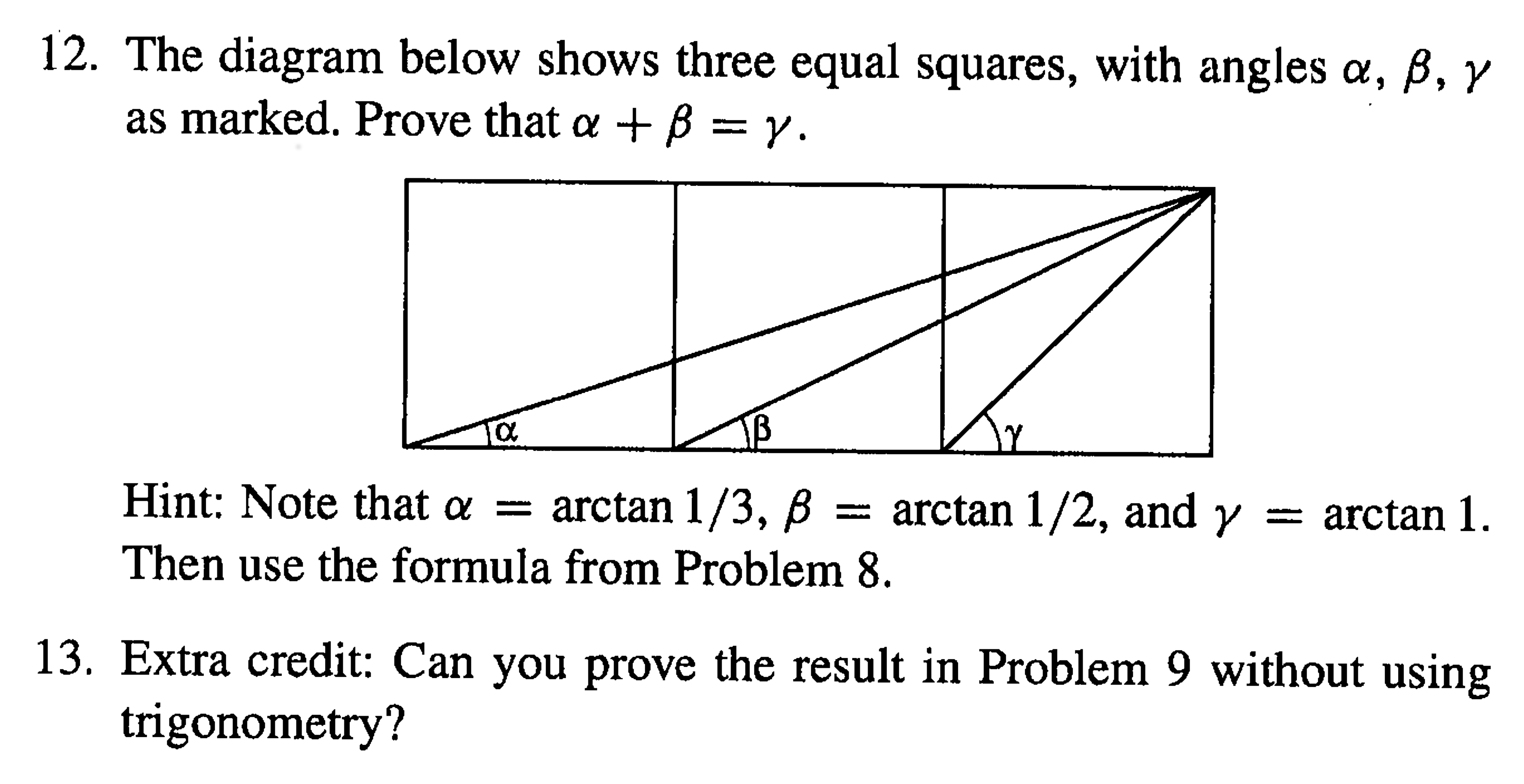}
\end{center}

Вот ещё несколько примеров из русской книжки.
\medskip

Объяснение тригонометрических функций произвольного угла: 

\begin{center}
\includegraphics[width=0.7\textwidth]{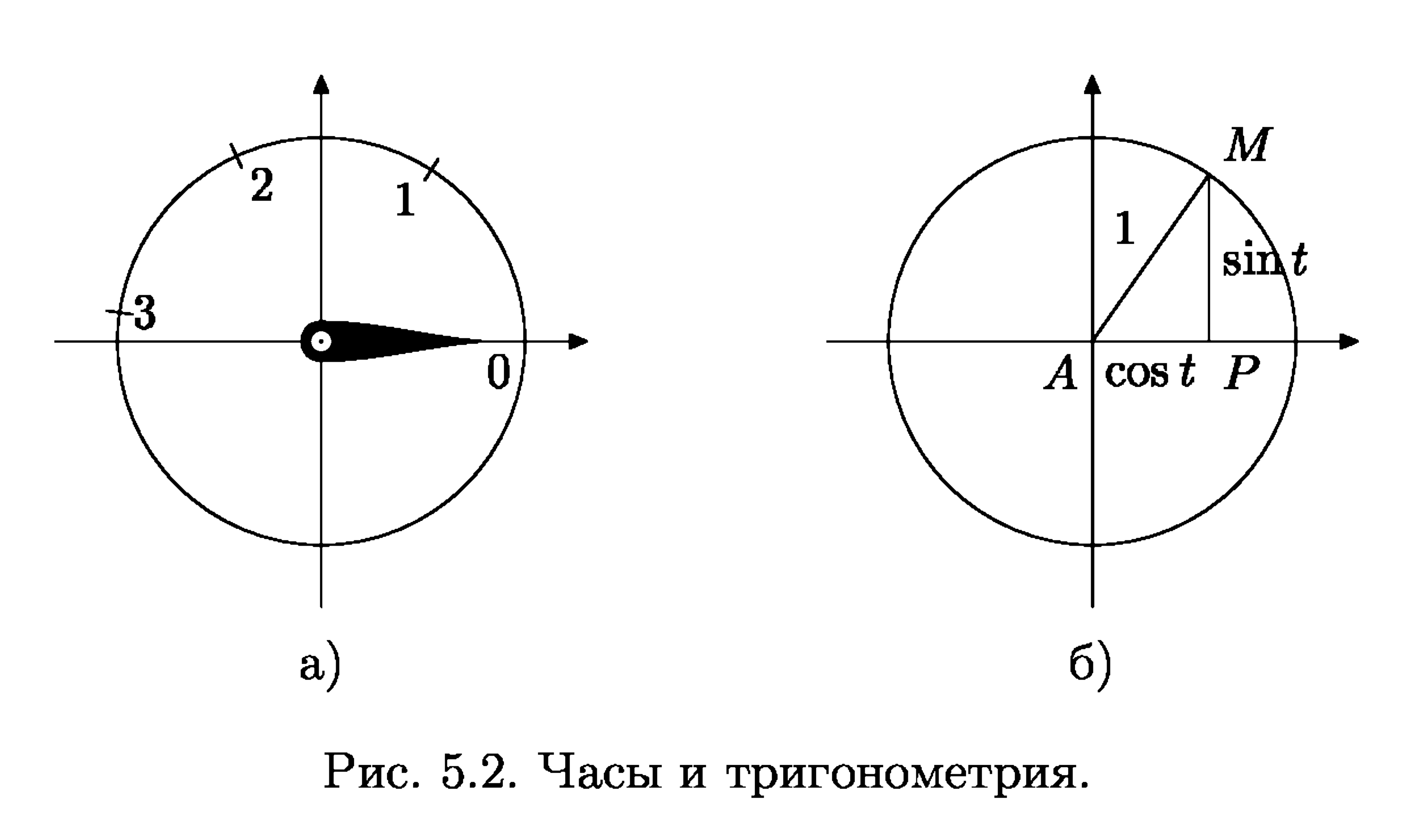}
\end{center}

Постановка вопроса, ответ на который получится дальше:

\begin{center}
\includegraphics[width=0.7\textwidth]{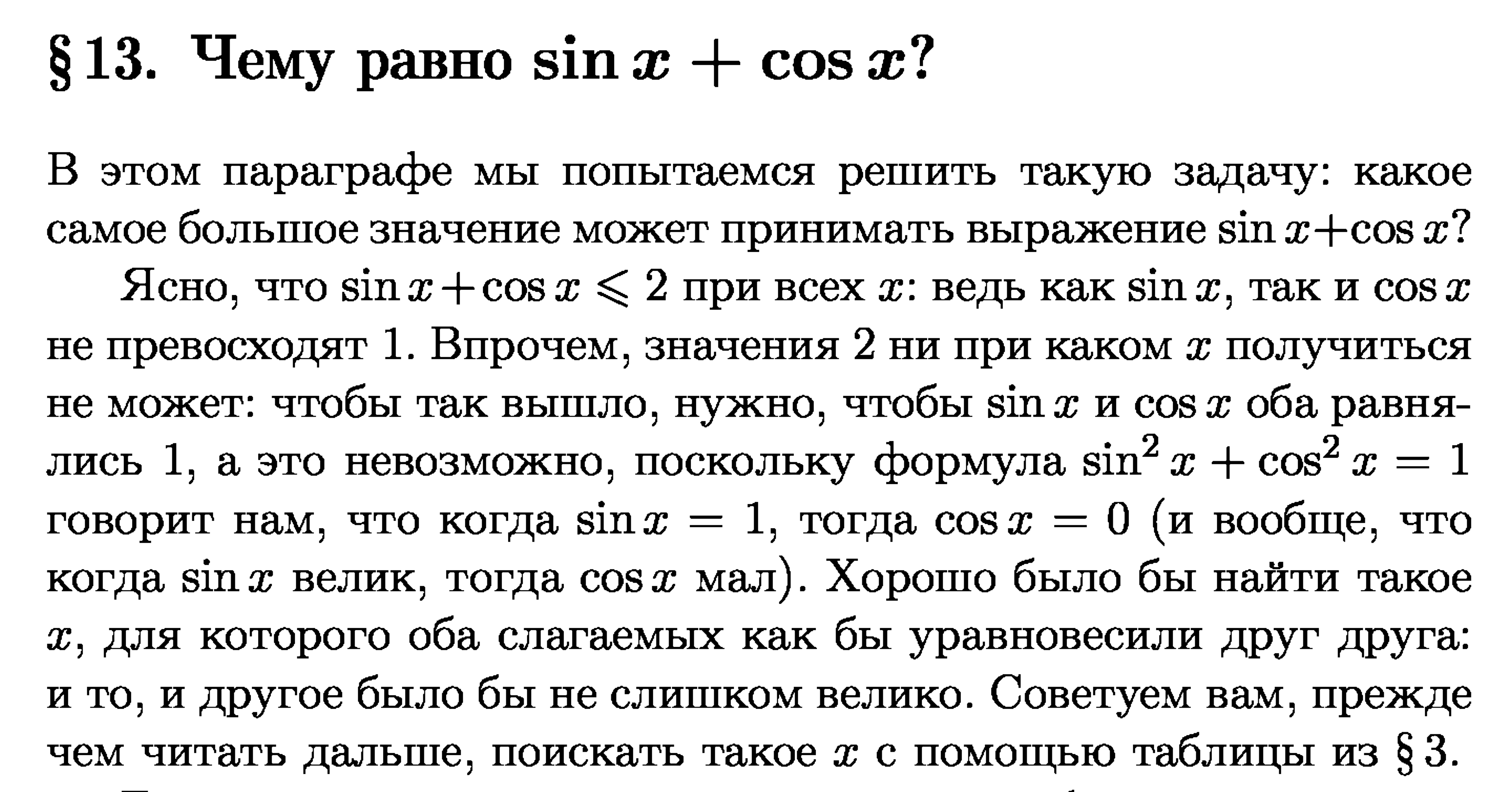}
\end{center}

\subsection*{Локальное качество текста}

Мы привели разные примеры задач из книжек ВЗМШ, но тут важен не только выбор тем и задач, но и, так сказать, <<локальное качество текста>>, о котором часто (я бы даже сказал: обычно) забывают, обсуждая книжки по математике для школьников.

Под словами <<локальное качество>> я имею в виду: насколько текст можно (или нельзя) улучшить, не меняя его содержания, только редактируя небольшие куски. И.\,М. не употреблял такого выражения, но объяснял мне, что важно не только то, что пишешь, но и то, как, и что можно писать хорошо, но по-разному (приводя в пример Пушкина и Толстого,\footnote{%
Кстати, о Толстом: в << Войне и мире>> есть фрагмент, объясняющий парадокс про Ахиллеса и черепаху, и написано там очень грамотно и качественно:
\begin{quote}
Известен так называемый софизм древних, состоящий в том, что Ахиллес никогда не догонит впереди идущую черепаху, несмотря на то, что Ахиллес идет в десять раз скорее черепахи: как только Ахиллес пройдет пространство, отделяющее его от черепахи, черепаха пройдет впереди его одну десятую этого пространства; Ахиллес пройдет эту десятую, черепаха пройдет одну сотую и т. д. до бесконечности. Задача эта представлялась древним неразрешимою. 
\end{quote}
Современные авторы (даже выпускники 57 школы, правда, <<гуманитарного класса>>, и вообще не первого сорта) далеко не так грамотны (цитата из ЖЖ \textbf{tema}):
\begin{quote}
Мой брательник (alexeilebedev) пишет: <<Может ли мне кто-нибудь объяснить, в чём проблема с определением правдивости утверждения `все критяне лжецы', сделанного обитателем Крита?>>.

Меня тоже раздражает этот псевдопарадокс, но он находится на месте номер пятьдесят по сравнению с другим, который меня выбешивает с детства --- про Ахиллеса и черепаху. Про то, что Ахиллес никогда не догонит черепаху (\url{http://ru.wikipedia.org/wiki/Ахиллес_и_черепаха}).

$\langle\ldots\rangle$, я не могу передать, как меня эта восточная тупость плющит. Если черепаха должна проползти половину своих десяти метров, то $\langle\ldots\rangle$ Ахиллес должен пробежать половину своей $\langle\ldots\rangle$  стометровки. И так до бесконечности.

Собиратель историй решил $\langle\ldots\rangle$ нас на базаре: это не апория Зенона, а $\langle\ldots\rangle$ $\langle\ldots\rangle$.
\end{quote}
(Опущена обсценная лексика, но от неё яснее не становится.)} как я уже говорил).
 
Оценка учебников по традиции игнорирует эту сторону дела. Даже комиссия Академии наук по оценке учебников по математике под руководством В.\,А.\,Васильева, который долго героически занимался этим неблагодарным делом, пока это имело хоть какой-то смысл, не могла этого делать:  по статусу комиссии её работа была ограничена математическим содержанием учебников, и уже из-за этого (а также по обычным советско-российским причинам) по большей части пропала впустую: в большинстве случаев не это в первую очередь делало учебник непригодным, а именно локальное качество.

Трудно объяснить, что имеется в виду под <<локальным качеством>>, не приводя примеры (того, что кажется хорошим или плохим).  Попытаемся посмотреть на изложение одного и того же рассуждения в разных учебниках геометрии. В качестве примера возьмём теорему о вписанном угле. Вот как она излагается в книжке А.\,В.\,Погорелова <<Элементарная геометрия. Планиметрия>> (Москва, Наука, 1969):
\begin{figure}[p]
\begin{center}
\includegraphics[width=0.8\textwidth]{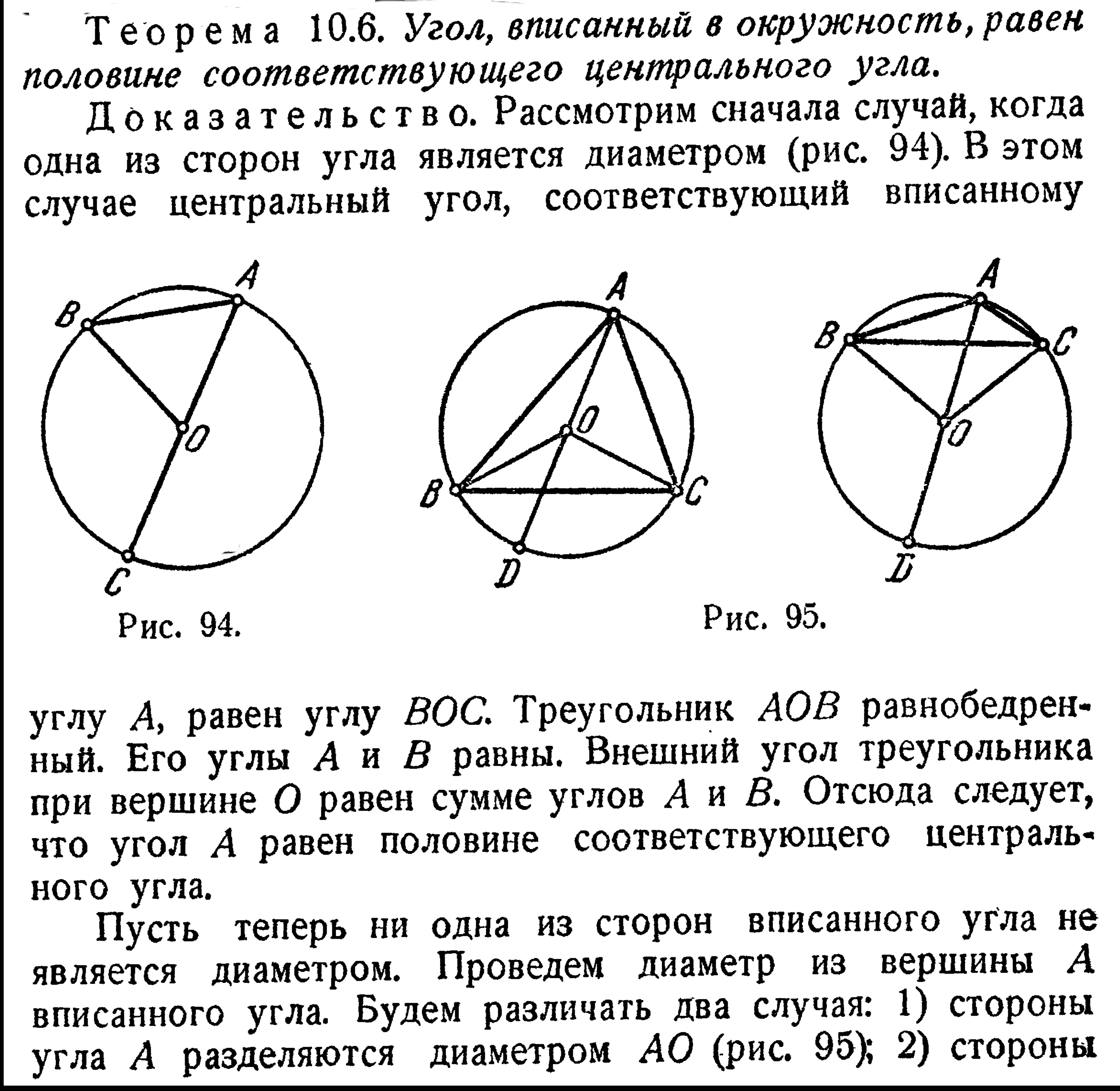}\\
\includegraphics[width=0.8\textwidth]{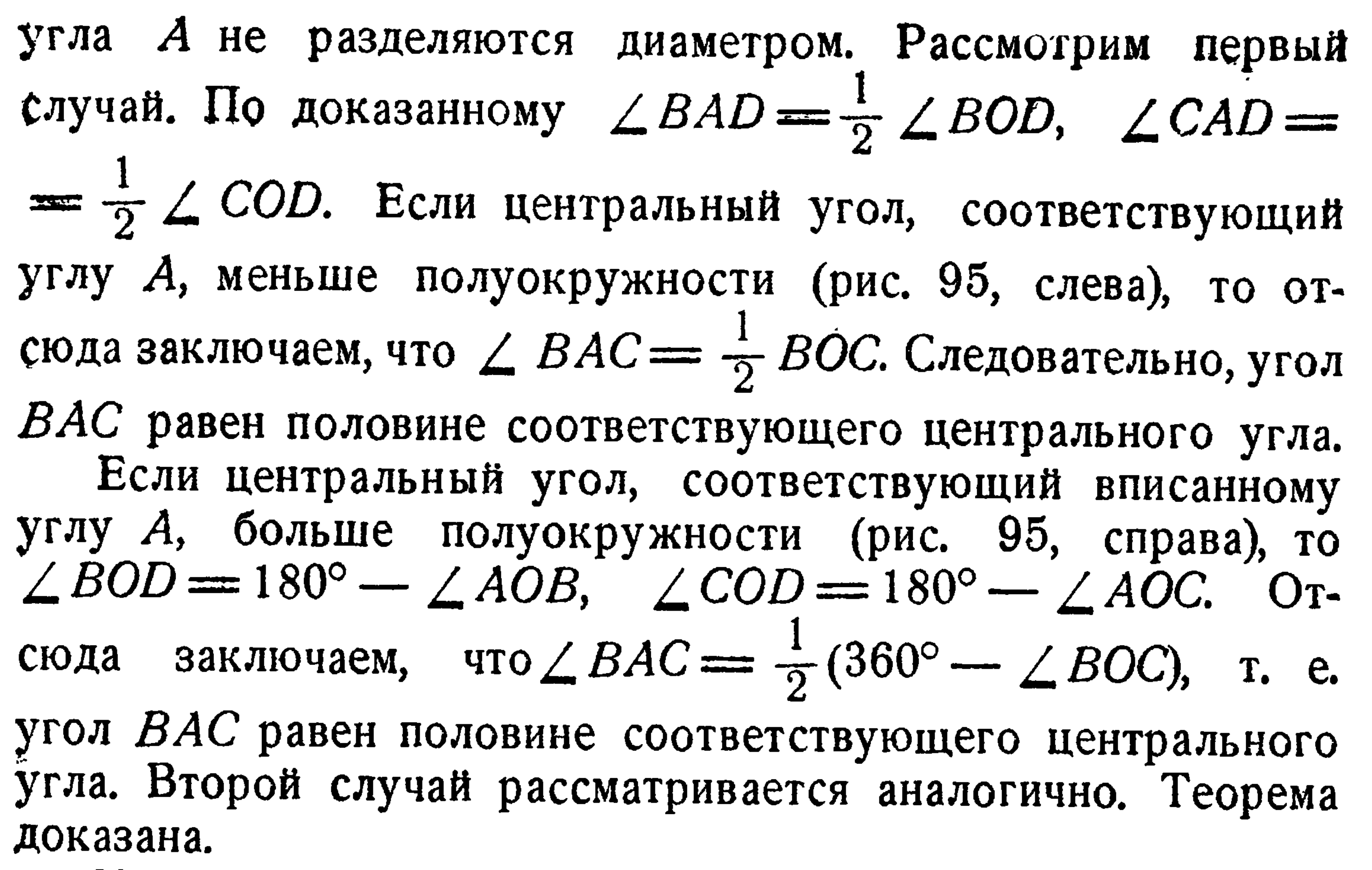}
\end{center}
\caption{Теорема о вписанном угле (<<Элементарная геометрия. Планиметрия>>, А.\,В.\,Погорелов).}\label{pogorelov-ugol}
\end{figure}
\clearpage

\begin{figure}[h]
\begin{center}
\includegraphics[width=0.67\textwidth]{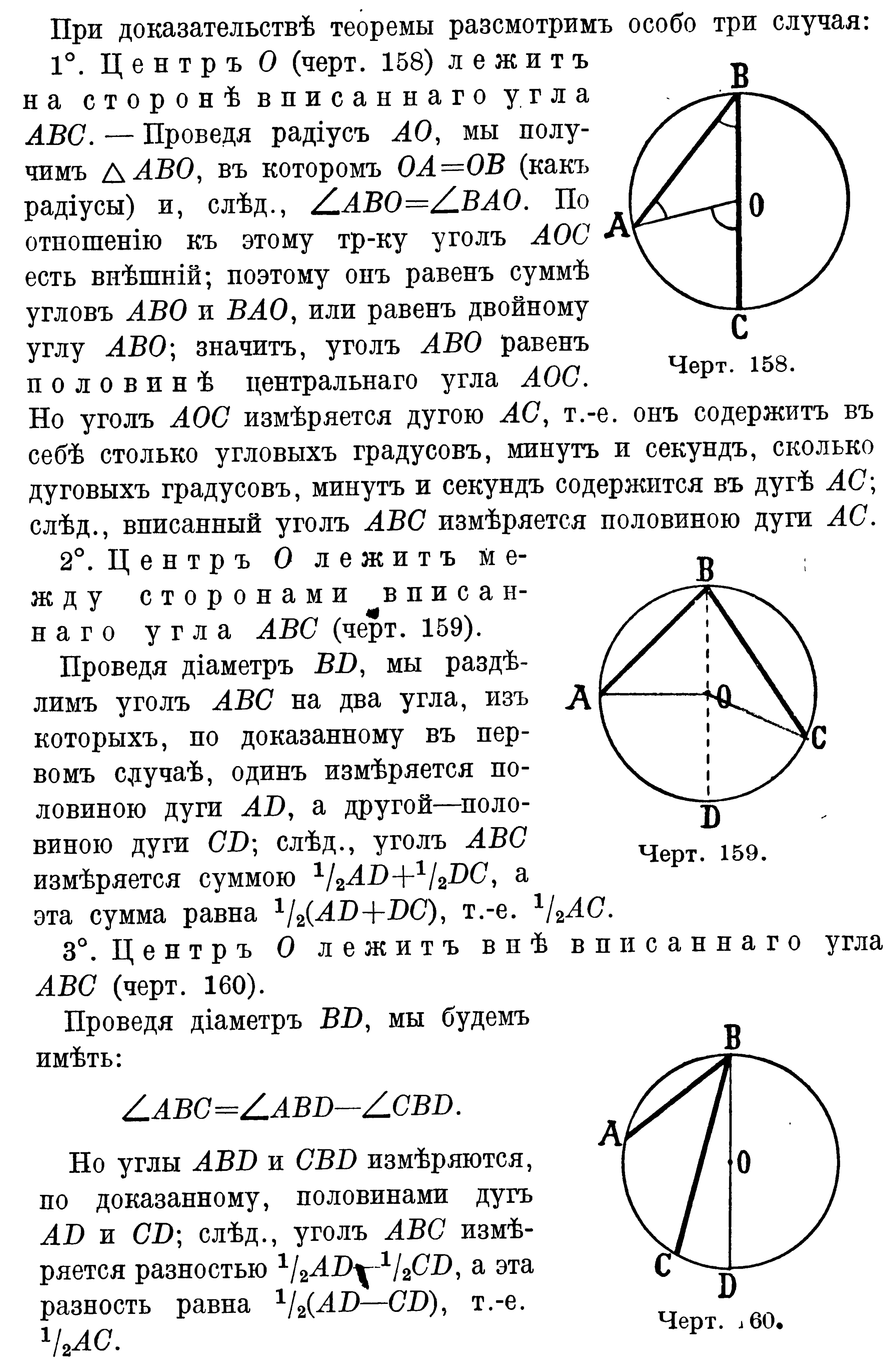}
\end{center}
\caption{Доказательство теоремы о вписанном угле в классическом учебнике Киселёва.}\label{kiselev-ugol}
\end{figure}

Посмотрим на то же самое доказательство у Киселёва (рис.~\ref{kiselev-ugol}).
\clearpage

На первый взгляд никакой разницы, но посмотрим подробнее. Для начала: у Киселёва рисунки расположены рядом с соответствующими рассуждениями. Затем: в первом же предложении говорится не просто <<одна из сторон угла является диаметром>> (и школьник должен понять, о каком угле на рисунке идёт речь), а конкретно <<центр $O$ лежит на стороне вписанного угла $ABC$>>. Равные углы при основании равнобедренного треугольника отмечены на рисунке. (Кстати, у Погорелова треугольник на рисунке выглядит равносторонним.) Вместо <<соответствующего центрального угла>> сказано конкретно <<центрального угла $AOC$>> и т.\,д.

Дальше у Погорелова идёт довольно запутанный текст, связанный с тем, что центральный угол, больший $180^\circ$, определяется как дополнение до $360^\circ$ дополнительного к нему угла (это связано с попытками навести строгость --- отдельная тема для обсуждения, в которое мы вдаваться не будем). Но и оно изложено плохо: в ключевом шаге (<<Отсюда заключаем>> в четвёртой строке снизу) пропущено главное: что угол $BAC$ равен сумме углов $BAO$ и $CAO$.

В конце (предпоследняя строка) у Погорелова написано <<Второй случай рассматривается аналогично>> --- но не так просто понять, что речь идёт о случае, когда <<стороны угла не разделяются диаметром>>, и что на картинке его искать бесполезно --- там показаны два случая, но другие. (К тому же и картинка, и перечисление случаев попали на предыдущую страницу.) Можно ещё заметить, что на рисунке Погорелова проведена совершенно лишняя линия $BC$.

Мы воспроизвели текст из книжки Погорелова, которая ещё не была учебником: при подготовке учебника две его книжки (про планиметрию и стереометрию) были объединены в одну и переработаны. Это тоже отдельная история; сейчас мы посмотрим только на то же самое доказательство в объединённой книжке (рис.~\ref{pogorelov2-ugol}):

\begin{figure}[p]
\begin{center}
\includegraphics[width=0.6\textwidth]{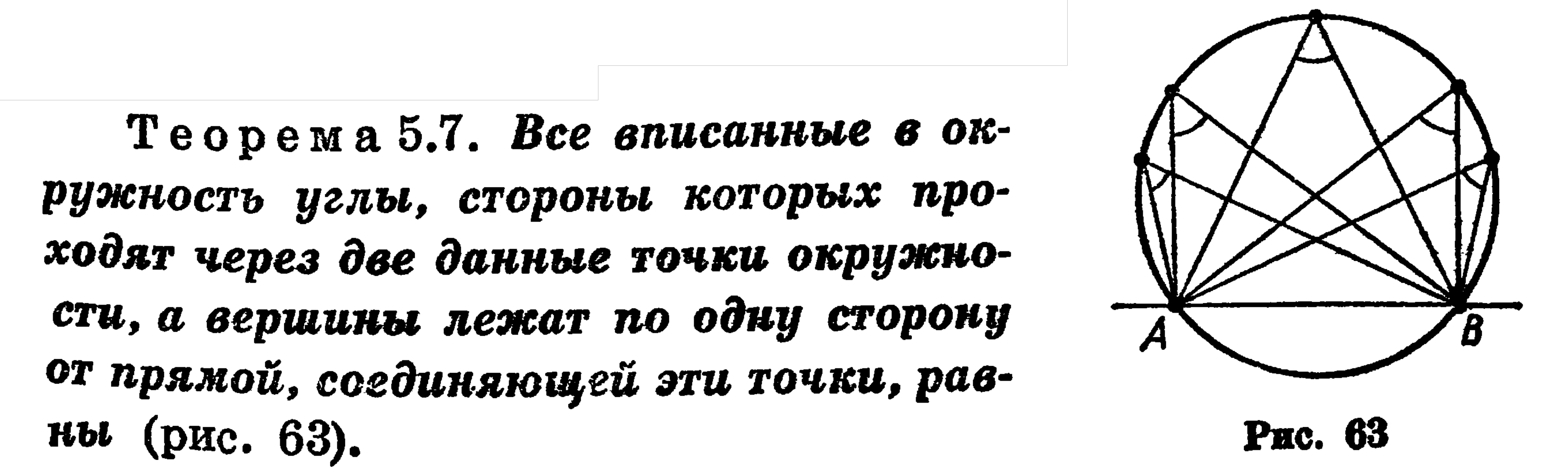}\\
\includegraphics[width=0.6\textwidth]{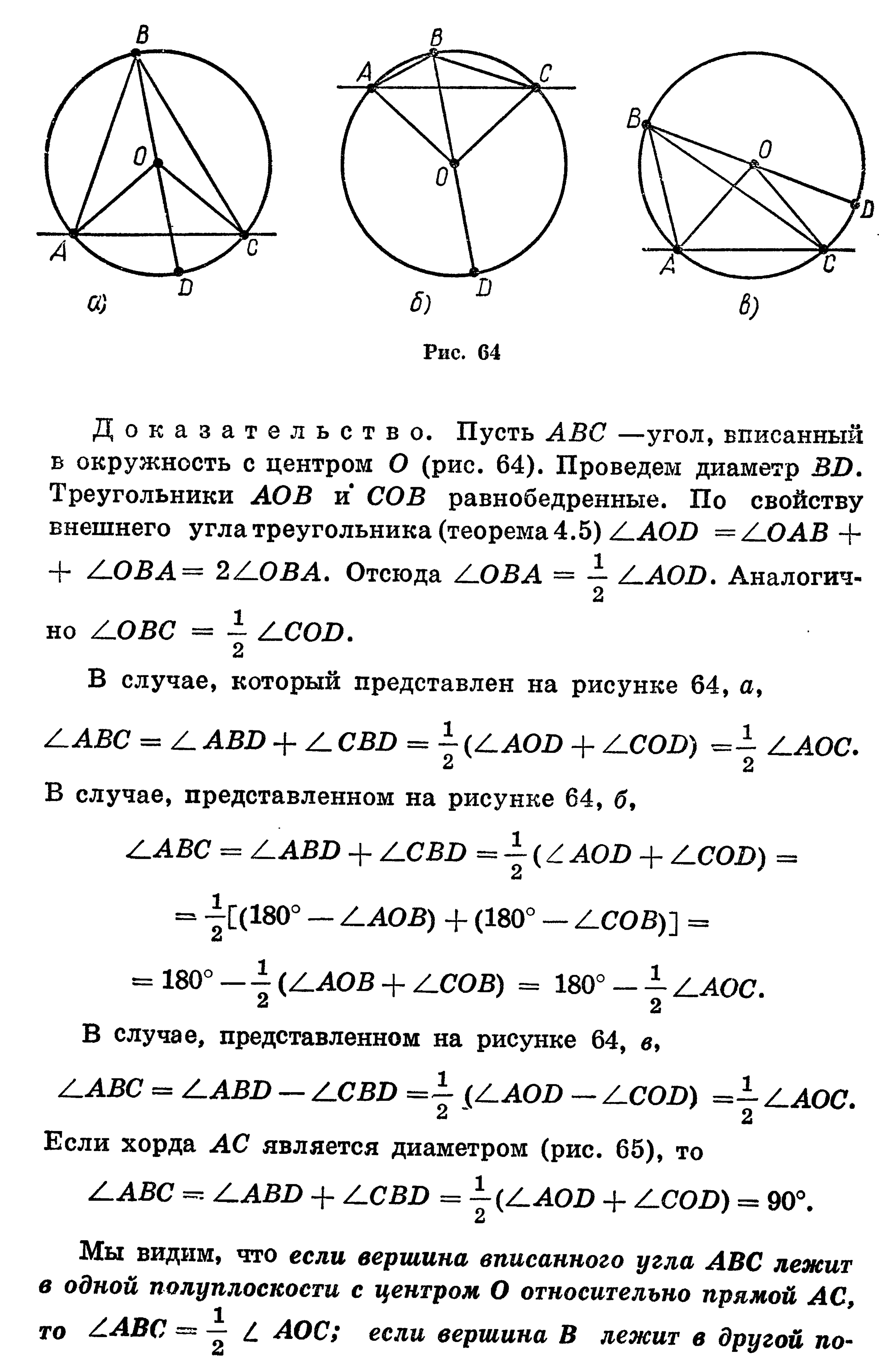}\\
\includegraphics[width=0.6\textwidth]{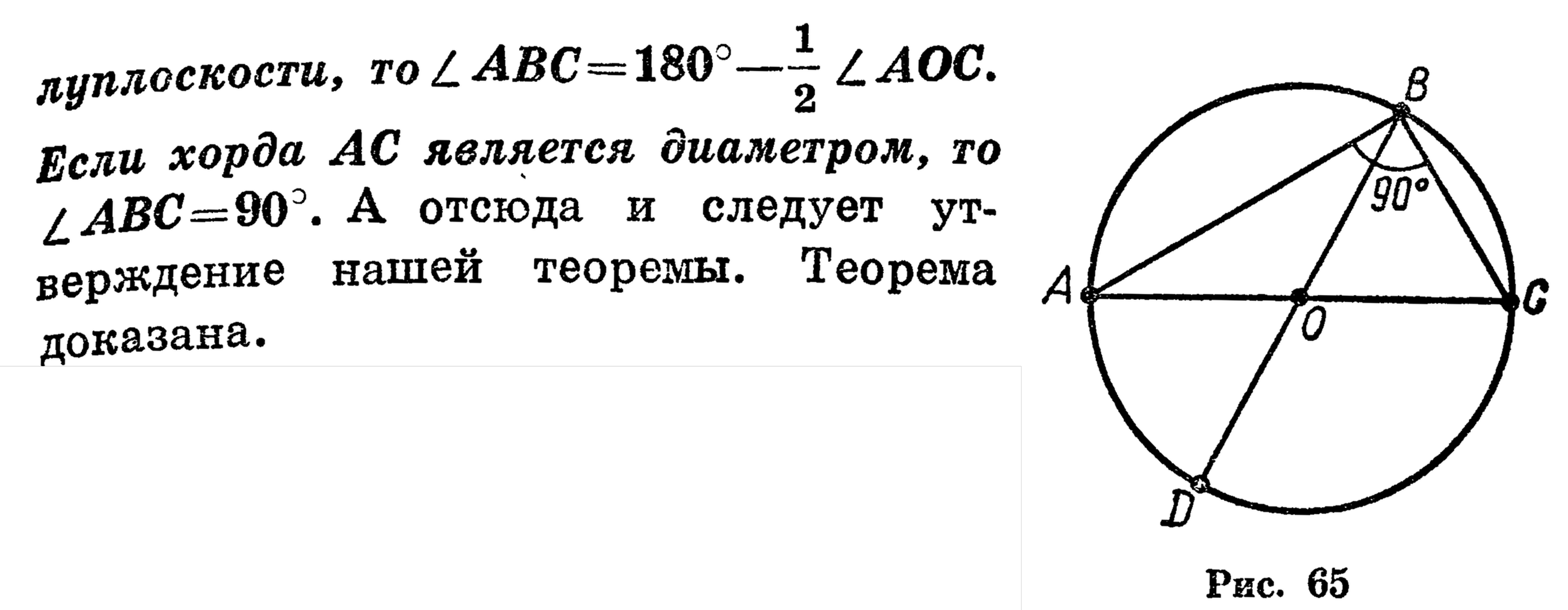}\\
\end{center}
\caption{Теорема о вписанном угле (<<Геометрия 6--10>>, А.\,В.\,Погорелов, Москва, Просвещение, 1981.)}
\label{pogorelov2-ugol}
\end{figure}
\clearpage

Прежде всего, тут изменена формулировка: говорится не об одном вписанном угле, который равен половине дуги, а о том, что все такие углы равны. Но доказательство начинается как раньше: <<пусть $ABC$ --- угол, вписанный в окружность>>, и если читатель вообще за чем-то следит в тексте, то у него сразу же возникнет вопрос: а где же здесь разные углы, о равенстве которых говорится в теореме?  Дальше первый случай (когда одна из сторон угла проходит через центр) не разбирается отдельно, а погружён внутрь двух других случаев. К тому же говорится, что <<треугольники $AOB$ и $COB$ равнобедренные>>, но никаких выводов о равенстве углов при основании не делается, а это подразумевается само собой, и то не сразу, а только после предложения о внешнем угле. Опять же на картинке не отмечены равные углы.  Напечатанный курсивом в конце текст вряд ли читатель соотнесёт с чем-то из сказанного ранее (те два случая, которые в нём упоминаются, никак не соотносятся с разобранными раньше случаями, а рисунок смешивает то и другое).%
\footnote{%
Надо сказать, что в более поздних изданиях (я видел издание 2009 года) текст сильно улучшен: случай с диаметром возвращён, лишняя линия на картинке убрана, осталось только углы пометить.}

\bigskip
Наконец, в качестве почти что эталонного плохого текста можно привести пример из книжки Л.\,С.\,Пон\-тря\-ги\-на <<Математический анализ для школьников>> (Москва, Наука, 1980). Понтрягин пишет в предисловии: <<Эта небольшая книга, объёмом около пяти листов, рассчитана на то, чтобы при удаче стать учебником математического анализа в средней школе>>. 

\begin{figure}[p]
\begin{center}
\includegraphics[width=0.8\textwidth]{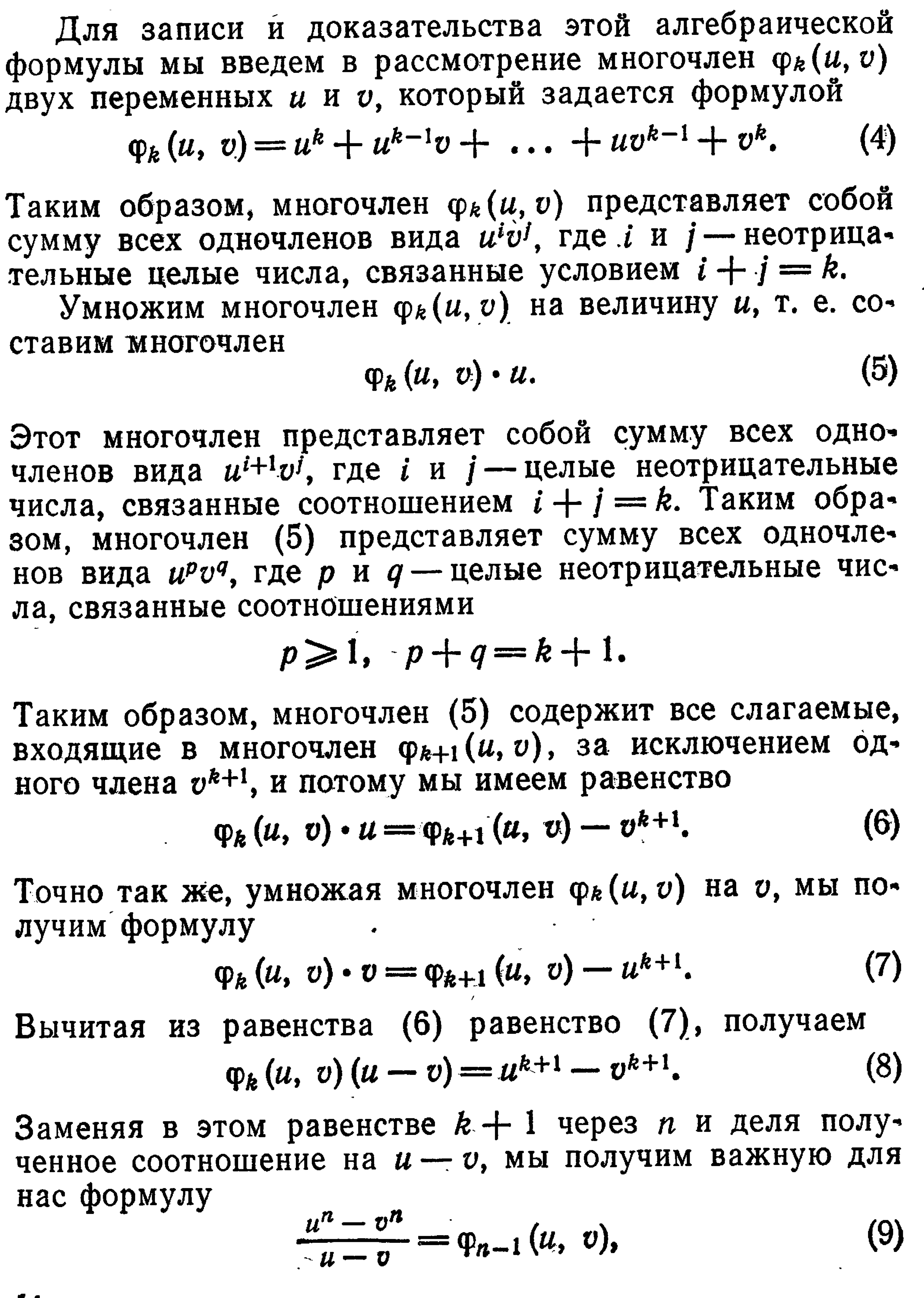}\\
\includegraphics[width=0.8\textwidth]{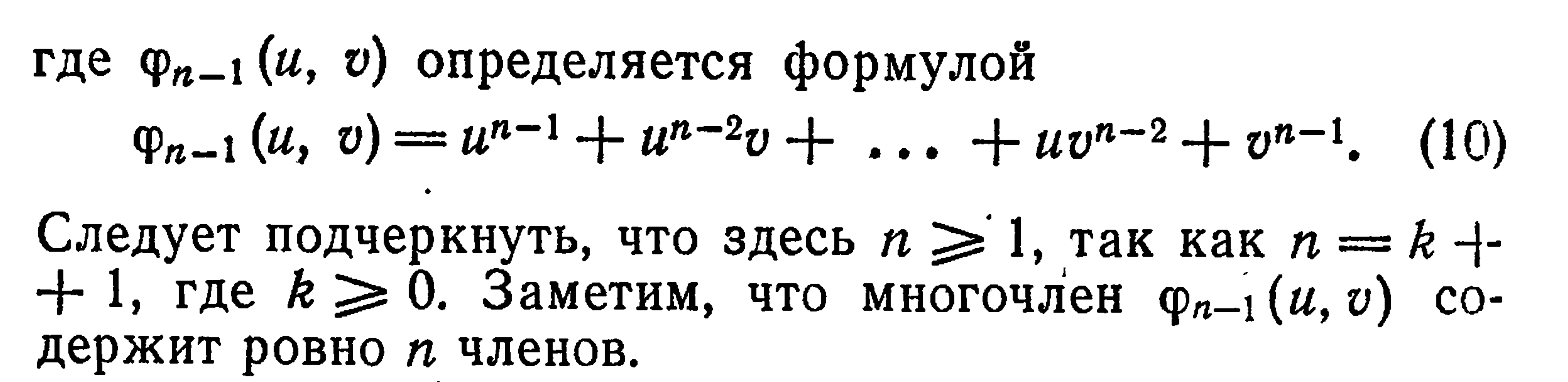}\\
\end{center}
\caption{Разложение $u^n-v^n$ на множители, изложение Л.\,С.\,Понтрягина.}
\label{pontryagin}
\end{figure}

В некоторый момент речь заходит о доказательстве тождеств вида 
\begin{align*}
u^2-v^2&=(u-v)(u+v)\\
u^3-v^3&=(u-v)(u^2+uv+v^2)\\
u^4-v^4&=(u-v)(u^3+u^2v+uv^2+v^3),\\
  &\ldots
\end{align*}
в общем виде (для произвольных степеней); они нужны для дифференцирования функции $x\mapsto x^n$. Вот как это доказательство излагается (рис.~\ref{pontryagin}).

\clearpage

<<Теперь пишут\ldots\ Необыкновенно пишут! Возьмешь. Раз прочтешь. Нет! Не понял. Другой раз --- то же. Так и отложишь в сторону\ldots>> (М.\,А.\,Булгаков)

\subsection*{Послесловие}

Разные куски этого текста были написаны в 2003 -- 2013 годах (в частности, для юбилейной конференции, посвящённой столетию со для рождения Гельфанда). С тех пор многое переменилось --- и, в частности, вклад Гельфанда в создание и поддержание высокого уровня советской математики, включая военно-прикладную, который казался тогда его безусловной заслугой, теперь уже вызывает более сложные чувства. Но моя личная благодарность И.\,М. за уроки и образцы в области преподавания и написания текстов для школьников и студентов остаётся неизменной.

Автор благодарен С.\,М.\,Львовскому и В.\,Шувалову, прочитавшим исходный вариант текста, за предложения и замечания.

\end{document}